\documentclass[11pt]{amsart}
\usepackage[margin=1.15in]{geometry}
\usepackage{amssymb}
\usepackage{xcolor}
\definecolor{linkblue}{HTML}{22506A}
\usepackage[colorlinks=true,linkcolor=linkblue,citecolor=linkblue,urlcolor=linkblue]{hyperref}
\usepackage{stmaryrd}
\usepackage{enumitem}
\usepackage{tikz-cd}
\usepackage{microtype}
\DeclareMathOperator{\im}{im}
\DeclareMathOperator{\Pic}{Pic}
\DeclareMathOperator{\CH}{CH}
\DeclareMathOperator{\id}{id}
\DeclareMathOperator{\Ann}{Ann}
\DeclareMathOperator{\GL}{GL}
\DeclareMathOperator{\PGL}{PGL}
\DeclareMathOperator{\Sym}{Sym}
\DeclareMathOperator{\Bl}{Bl}
\DeclareMathOperator{\Std}{Std}
\DeclareMathOperator{\Nm}{Nm}
\DeclareMathOperator{\ddiv}{div}
\DeclareMathOperator{\sgn}{sgn}
\DeclareMathOperator{\Hilb}{Hilb}
\DeclareMathOperator{\Aut}{Aut}
\DeclareMathOperator{\Gr}{Gr}
\DeclareMathOperator{\rig}{rig}

\newcommand{\mg}[1]{\mathcal{M}_{#1}}
\newcommand{\mgbar}[1]{\overline{\mathcal{M}}_{#1}}
\newcommand{\rigidify}{\mkern-5mu\fatslash\mkern5mu}

\newtheorem{theorem}{Theorem}[section]
\newtheorem{maintheorem}{Theorem}

\newtheorem{lemma}[theorem]{Lemma}
\newtheorem{proposition}[theorem]{Proposition}
\newtheorem{corollary}[theorem]{Corollary}

\theoremstyle{definition}
\newtheorem{definition}[theorem]{Definition}
\newtheorem{note}[theorem]{Note}

\numberwithin{equation}{section}

\title{The integral Chow rings of $\mathcal A_3$ and $\mathcal M_3$}
\author{Emma Bishop}

\begin{document}

\begin{abstract}
Over $\mathbb C$, we compute the integral Chow ring of $\mathcal A_3$, the moduli stack of
principally polarized abelian threefolds. Along the way, we compute the integral Chow ring of
$\mathcal M_3$, the moduli stack of smooth genus-three curves. The ring $\operatorname{CH}(\mathcal A_3)$
is generated by the Chern classes of the Hodge bundle, along with three classes supported on the
locus of decomposable abelian varieties, subject to 16 relations. The ring $\operatorname{CH}(\mathcal M_3)$
is generated by the Chern classes of the Hodge bundle, with five relations. In particular,
we find that the top Chern class of the Hodge bundle, $\lambda_3$, has order 504 on $\mathcal A_3$
and 252 on $\mathcal M_3$. For $\mathcal A_3$, this removes the factor of 4 from the bound
$252\mid\operatorname{ord}(\lambda_3)\mid 1008$ of Ekedahl and van der Geer.
\end{abstract}

\maketitle

\section{Introduction}

Since the introduction of equivariant intersection theory by Totaro \cite{Totaro} and
Edidin and Graham \cite{EG98}, integral Chow rings have been computed for several
important moduli spaces, including $\mg2$ \cite{Vistoli}, $\mgbar2$ \cite{Larson},
the stacks $\mathcal H_g$ of hyperelliptic curves (\cite{DL} for odd genus, \cite{EF09} for even genus),
the stack $\mg3\setminus\mathcal H_3$ of smooth non-hyperelliptic curves of genus
three \cite{DLFV}, and $\mg2^{\mathrm{ct}}$ \cite{HJL}. Additionally, the ``almost integral''
Chow ring, with $\mathbb Z[1/6]$ coefficients, of $\mgbar3$ was computed in \cite{Pernice}.
For moduli spaces of principally polarized abelian varieties,
much less is known about their integral Chow rings: $\mathcal A_1$, using
$\mathcal A_1\cong\mg{1,1}$, has been known since \cite{EG98}, $\mathcal A_2$, using
$\mathcal A_2\cong\mg2^{\mathrm{ct}}$, is due to \cite{HJL}, and $\mathcal A_3$ is known rationally
due to \cite{vdG}. The order of $\lambda_3$ on $\mathcal A_3$ is known up to a factor of 4 due to
\cite{EvdG}, who showed that $252\mid\operatorname{ord}(\lambda_3)\mid 1008$.

In this paper we compute the integral Chow ring of $\mathcal A_3$, and as a
consequence of our methods we also determine the integral Chow ring of $\mg3$.
In particular, the order of $\lambda_3$ is 504 on $\mathcal A_3$ and 252 on $\mg3$.

\subsection{Main results}

For both $\mg3$ and $\mathcal A_3$, let $\mathbb E_3$ denote the Hodge bundle and put
$$
\lambda_1=c_1(\mathbb E_3),\qquad \lambda_2=c_2(\mathbb E_3),
\qquad \lambda_3=c_3(\mathbb E_3).
$$
Let $\mathcal A_3^{\mathrm{ind}}$ denote the indecomposable locus of $\mathcal A_3$.
Let $A_{1,1,1}\subseteq\mathcal A_3$ be the closed substack of products of three elliptic curves,
with inclusion $p_{111}:A_{1,1,1}\rightarrow\mathcal A_3$. Note that
$A_{1,1,1}=\left[\mathcal A_1^3/S_3\right]$. Let
$p_{12}:\mathcal A_1\times\mathcal A_2\rightarrow\mathcal A_3$ be the product morphism. Set
$$
\alpha_{12}=(p_{12})_*1,\qquad \alpha_{12}^{(1)}=(p_{12})_*t,\qquad \alpha_{111}=(p_{111})_*1,
$$
where $t$ is the first Chern class of the Hodge bundle of $\mathcal A_1$.

\begin{maintheorem}
We have
$$
\CH(\mathcal A_3)=\frac{\mathbb Z[\lambda_1,\lambda_2,\lambda_3,\alpha_{12},\alpha_{12}^{(1)},
\alpha_{111}]}{I_{16}},
$$
where
$\deg(\lambda_1,\lambda_2,\lambda_3,\alpha_{12},\alpha_{12}^{(1)},\alpha_{111})=(1,2,3,2,3,3)$
and the relations in $I_{16}$ are listed in Appendix \ref{A_3 relations appendix}. In particular,
the order of $\lambda_3$ is 504.
\end{maintheorem}

\begin{maintheorem}
We have
$$
\CH(\mg3)=\frac{\mathbb Z[\lambda_1,\lambda_2,\lambda_3]}
{(252\lambda_1^2, 24\lambda_2-12\lambda_1^2, 8\lambda_1\lambda_2+80\lambda_1^3,
28\lambda_3+28\lambda_1^3, 2\lambda_1\lambda_3+2\lambda_1^2\lambda_2+4\lambda_1^4)}.
$$
Pullback along the Torelli morphism is surjective and identifies
$$
\CH(\mg3)=\CH(\mathcal A_3)/(\alpha_{12}, \alpha_{12}^{(1)}, \alpha_{111},252\lambda_1^2)=
\CH(\mathcal A_3^{\mathrm{ind}})/(252\lambda_1^2).
$$
In particular, the order of $\lambda_3$ is 252.
\end{maintheorem}

\subsection{Outline of the proof}
The open locus $\mathcal A_3^{\mathrm{ind}}$ consists of abelian threefolds which are Jacobians
of smooth genus-three curves. The moduli stack $\mg3$ of smooth genus-three curves
contains a closed substack, $\mathcal H_3$, of hyperelliptic curves. This naturally creates
a corresponding closed substack in $\mathcal A_3^{\mathrm{ind}}$ of threefolds which
are Jacobians of curves from $\mathcal H_3$. Call this substack $\mathcal A_3^{\mathrm{hyp}}$,
and set
$\mathcal A_3^{\mathrm{nh}}=\mathcal A_3^{\mathrm{ind}}\setminus\mathcal A_3^{\mathrm{hyp}}$,
the non-hyperelliptic locus.
Let $A_{1,2}^\circ\subseteq\mathcal A_3$ be the locus of products of an elliptic curve and the Jacobian
of a smooth genus-two curve. Then
$$
\mathcal A_3=\mathcal A_3^{\mathrm{nh}}\sqcup\mathcal A_3^{\mathrm{hyp}}
\sqcup A_{1,2}^\circ\sqcup A_{1,1,1}.
$$

\textit{Step 1} (Section \ref{plane quartics section}).
Let $P$ be the space of plane quartics and $Z\subseteq P$ the locus of double conics.
We first deform $P$ to the normal cone of $Z$, then remove the part of the special fiber over singular
conics, and take the quotient $X$ by $\PGL_3\times\mathbb G_m$. The special divisor of $X$ is a
vector bundle over $B\PGL_2\times B\mathbb G_m$ and the top Chern class of its normal bundle
isn't a zero-divisor, so by the patching lemma (Lemma \ref{patching lemma})
$\CH(X)$ is a fiber product of the Chow rings
of the two strata (Proposition \ref{first patch}).

\textit{Step 2} (Sections \ref{envelope section}--\ref{rigidification Chow ring}).
Let $J\subseteq X$ be the open substack parametrizing smooth quartics and squarefree binary
octics on smooth conics. We compute $\CH(J)$ by excising the complement $B=X\setminus J$.
A partial analysis of $B$ produces six relations, which may or may not generate all relations.
Rather than try to continue the excision, we prove that every cycle pushed forward from $B$
is generated by these six relations after potentially multiplying by 3 (Theorem \ref{degree three}).
This tells us that the obstruction to our list being complete is annihilated by 3.
Restriction to the locus of smooth quartics shows that it is also annihilated by 112 in degrees
at least two (Lemma \ref{112 annihilator}). Since 3 and 112 are coprime,
the six relations are complete (Proposition \ref{excision result}).
We also determine the kernel of restriction to the divisor $H$ (Proposition \ref{H kernel}).

\textit{Step 3} (Section \ref{root section}).
We show that the stack $J$ is the rigidification of $\mathcal A_3^{\mathrm{ind}}$
by the involution (Proposition \ref{A_3 rigidification}). In fact, $\mathcal A_3^{\mathrm{ind}}$ is the
root gerbe of $L=\mathcal O_J(H)$, where $H\subseteq J$ is the locus of hyperelliptic curves.
We could have avoided taking a root gerbe if, in the deformation to the normal cone, we
instead had $\mathbb G_m$ act with weight 2, but proceeding the way that we have creates the
following serendipity: $\mg3$ is the root \textit{stack} of $J$ along $H$ (Propositions
\ref{A_3^ind root gerbe} and \ref{M_3 root stack}). Therefore we may compute the integral Chow
rings of both $\mathcal A_3^{\mathrm{ind}}$ and $\mg3$ by just applying the root gerbe and root
stack formulas.

\textit{Step 4} (Section \ref{A_12 Chow section}).
We next compute the integral Chow ring of $A_{1,2}^\circ=\mg{1,1}\times\mg2$ and its normal
bundle in $\mathcal A_3\setminus A_{1,1,1}$. Unfortunately, the top Chern class of this bundle
is a zero-divisor. Fortunately, its annihilator is very well-behaved: it is just $120\mathbb Z$,
concentrated in degree zero (Theorem \ref{A_12 annihilator}).

\textit{Step 5} (Sections \ref{BG_3 section}--\ref{A_1,1,1 normal bundle section}).
Let $G_3=\mathbb G_m^3\rtimes S_3$.
We take inspiration from \cite[\S 5]{Larson} to compute $\CH(BG_3)$, and then use
excision to compute $\CH(A_{1,1,1})$ (Propositions \ref{BG Chow ring} and \ref{A_1,1,1 Chow ring}),
which we do by exhibiting $A_{1,1,1}$ as a quotient of an open in $(\mathbb A^2)^3$ by $G_3$.
The top Chern class of the normal bundle of $A_{1,1,1}$ is again a zero-divisor.

\textit{Step 6} (Sections \ref{M_2^ct test family section}--\ref{final ring section}).
We use $\mg2^{\mathrm{ct}}$ as a test family, and show that the simultaneous restriction of
$\CH(\mathcal A_3)$ to each of $\CH(\mathcal A_3^{\mathrm{ind}})$, $\CH(A_{1,2}^\circ)$,
$\CH(A_{1,1,1})$, $\CH(\mg2^{\mathrm{ct}})$, and
$\CH(\mathcal A_3)_{\mathbb Q}$ is injective. After finding a set of generators for
$\CH(\mathcal A_3)$, this allows us to determine the relations and hence find a presentation
of the Chow ring of $\mathcal A_3$.

\subsection{Computational supplement}
There are many quite in-depth algebra computations which occur throughout the paper.
There is a computational supplement available in the ancillary files which
carries out these computations in Macaulay2.

\subsection{Conventions}
We work over $\mathbb C$. All Chow rings are integral unless otherwise stated. We will
occasionally conflate $\mg{1,1}$ with $\mathcal A_1$. For a vector bundle $E$ on a space
$X$, we take $\mathbb P(E)$ to be the space of lines in $E$. We will frequently use
$e$ to refer to the top Chern class of a vector bundle (its Euler class).

\subsection{AI disclosure}
The overall structure and mathematical ideas of the paper are the author's, and GPT-6 Astra assisted
with the ``dirty work'' in the excision computations. In particular, GPT-6 Astra wrote the entire
computational supplement, including all of the Macaulay2 code used to find presentations of rings,
compute kernels, find annihilators, etc. Great care was put into the exposition of this paper and
100\% of the writing was done by the author. Proofreading, editing, and typesetting of equations
and tables was assisted by Opus 5.5, GPT-6 Astra, and GPT-6.1 Sol.

\subsection{Acknowledgments}
The author benefited greatly from conversations with Veronica Arena, Giovanni Inchiostro, and
Michael Ruofan Zeng. The author would like to thank Yuchen Liu for suggesting this problem. The
idea of using the deformation to the normal cone originated from a conversation with Jarod Alper.

\section{Background}

\subsection{Background on $\mathcal A_g$}

We first state some standard results which will help us identify normal and excess bundles.
\begin{lemma}[Kodaira--Spencer, {\cite[Chapter III, \S9]{FC}}]\label{Kodaira-Spencer}
We have
$$
T_{\mathcal A_g}=\Sym^2(\mathbb E_g^\vee).
$$
\end{lemma}

\begin{lemma}\label{Kodaira-Spencer normal}
Let $g_1,\dots,g_r\geq1$, put $g=\sum_i g_i$, and let
$$
p:\prod_{i=1}^r\mathcal A_{g_i}\rightarrow\mathcal A_g
$$
be $(A_1,\dots,A_r)\mapsto A_1\times\cdots\times A_r$. Then the normal bundle along $p$ is
$$
N_p=\bigoplus_{i<j}\left(\mathbb E_{g_i}^\vee\otimes\mathbb E_{g_j}^\vee\right),
$$
where each Hodge bundle is pulled back from its corresponding factor.
\end{lemma}

\begin{proof}
Since $p^*\mathbb E_g=\bigoplus_i\mathbb E_{g_i}$,
Lemma \ref{Kodaira-Spencer} gives
$$
p^*T_{\mathcal A_g}=\bigoplus_i\Sym^2(\mathbb E_{g_i}^\vee)\oplus
\bigoplus_{i<j}\left(\mathbb E_{g_i}^\vee\otimes\mathbb E_{g_j}^\vee\right).
$$
The normal bundle is then the quotient of this by
$$
\bigoplus_{i=1}^rT_{\mathcal A_{g_i}}=\bigoplus_{i=1}^r\Sym^2(\mathbb E_{g_i}^\vee),
$$
which proves the lemma.
\end{proof}

\begin{definition}
The \textit{Torelli morphism} $\tau:\mg g\rightarrow\mathcal A_g$ is the morphism which sends
a smooth genus $g$ curve $C$ to its Jacobian, $J(C)$.
\end{definition}

We now define the important strata of $\mathcal A_3$.

\begin{definition}\label{strata definition}
Put
$$
A_{1,1,1}=[\mg{1,1}^3/S_3],\qquad A_{1,2}^\circ=\mg{1,1}\times\mg2,
$$
where both $A_{1,1,1}$ and $A_{1,2}^\circ$ map to $\mathcal A_3$ by taking products of
Jacobians. The hyperelliptic substack
$\mathcal A_3^{\mathrm{hyp}}\subseteq\mathcal A_3^{\mathrm{ind}}$
is the image of the stack $\mathcal H_3$ of smooth hyperelliptic curves of genus three under
the Torelli morphism. Set
$\mathcal A_3^{\mathrm{nh}}=\mathcal A_3^{\mathrm{ind}}\setminus\mathcal A_3^{\mathrm{hyp}}$,
the non-hyperelliptic locus. With this notation,
$$
\mathcal A_3=\mathcal A_3^{\mathrm{nh}}\sqcup\mathcal A_3^{\mathrm{hyp}}
\sqcup A_{1,2}^\circ\sqcup A_{1,1,1},
$$
and $\mathcal A_3^{\mathrm{hyp}}$, $A_{1,2}^\circ$, $A_{1,1,1}$ have codimension $1,2,3$ in
$\mathcal A_3$, respectively.
\end{definition}

\subsection{Intersection theory}
For a closed $Z$ with complement $U$ in $X$, the following lemmas will help us use the Chow
rings of $U$ and $Z$ to compute $\CH(X)$.

\begin{lemma}[Patching lemma]\label{patching lemma}
Let $X$ be a smooth variety endowed with the action of a group $G$,
let $i:Z\rightarrow X$ be a regular embedding of a smooth, closed, and $G$-invariant
subvariety, with normal bundle $N$.
Suppose that $c^G_{\mathrm{top}}(N)$ is not a zero-divisor in $\CH_G(Z)$, and let
$j:X\setminus Z\rightarrow X$ denote the open inclusion.
Let $R$ be the following Cartesian product of rings:
$$
\begin{tikzcd}
R\arrow[r, "i^*"]\arrow[d, "j^*"] & \CH_G(Z)\arrow[d, "q"]\\
\CH_G(X\setminus Z)\arrow[r, "p"] & \CH_G(Z)/(c^G_{\mathrm{top}}(N))
\end{tikzcd}
$$
where the bottom horizontal arrow $p$ is the \textit{specialization map} which sends the class of a
variety $V$ to the equivalence class $i^*z$ where $z$ is any element in the set
$(j^*)^{-1}([V])$.

Then there is an isomorphism
$$
\CH_G(X)\cong R.
$$
\end{lemma}

\begin{proof}
This is \cite[Lemma 3.4]{DLV21}.
\end{proof}

\begin{lemma}\label{fiber product lemma}
Let $R$ be a commutative ring and let $I,J\subseteq R$ be ideals. Then the natural map gives an
isomorphism
$$
\varphi:R/(I\cap J)\xrightarrow{\sim} R/I\times_{R/(I+J)}R/J.
$$
Moreover, if $L\subseteq J$ is an ideal and $x\in R$ is not a zero-divisor in $R/J$, then
$$
(L+(x))\cap J=L+xJ.
$$
\end{lemma}

\begin{proof}
The map $R\rightarrow R/I\times_{R/(I+J)}R/J$ has kernel $I\cap J$.
For surjectivity, let $(f\bmod I,g\bmod J)$ be an element of the fiber product. Then
$f-g=i+j$ for some $i\in I$ and $j\in J$, and $\varphi(f-i)=\varphi(g+j)=(f\bmod I,g\bmod J)$.
Hence $\varphi$ is an isomorphism.

For the second claim, suppose $l+xg\in J$ with $l\in L$. Since $L\subseteq J$, we have
$xg\in J$. As $x$ is not a zero-divisor in $R/J$, it follows that $g\in J$.
Thus $(L+(x))\cap J\subseteq L+xJ$, and the reverse inclusion is immediate.
\end{proof}

\begin{lemma}\label{Chow generation lemma}
Let $p:Z\rightarrow X$ be a closed immersion of smooth quotient stacks with complement
$j:U\rightarrow X$. Let $R$ be a ring with a graded homomorphism
$$
R\rightarrow\CH(X),
$$
and let $M\subseteq\CH(X)$ be a graded $R$-submodule. Suppose that $p_*1$ is in the image of
$R$, and that restriction induces surjections
$$
j^*|_M:M\twoheadrightarrow\CH(U), \qquad p^*|_M:M\twoheadrightarrow\CH(Z).
$$
Then $M=\CH(X)$.
\end{lemma}

\begin{proof}
Pick $\delta\in R$ whose image is $p_*1$. Since $p^*M=\CH(Z)$, the projection formula gives
$$
p_*\CH(Z)=p_*p^*M=(p_*1)M=\delta M\subseteq M.
$$
Let $\beta\in\CH(X)$. Then, since $j^*|_M$ is surjective, there exists $m\in M$ such that
$\beta=m+p_*z$ for some $z\in\CH(Z)$. Hence
$$
\CH(X)=M+p_*\CH(Z)=M+\delta M=M.
$$
\end{proof}

To aid in our patching, we will also use the following lemma, which is really more of an observation.
For our application, $\CH(\mathcal A_3)_{\mathbb Q}$ will play the role of $S$.

\begin{lemma}\label{test family restriction lemma}
Suppose that $X$ stratifies as $X=X_0\sqcup X_1\sqcup\dots\sqcup X_n$. Let $T$ be another space
with a map $T\rightarrow X$, and let $S$ be another ring with a map $\CH(X)\rightarrow S$.
Assume that all necessary pullbacks exist. If $\{x_{\alpha}\}$ is a generating set for $\CH(X)$ and
the simultaneous restriction
$$
\CH(X)\rightarrow\CH(X_0)\times\cdots\times\CH(X_n)\times\CH(T)\times S
$$
is injective, then $\CH(X)=\mathbb Z[x_{\alpha}]/K$, where $K$ is the kernel of the composition
$$
\mathbb Z[x_{\alpha}]\rightarrow\CH(X)\rightarrow
\CH(X_0)\times\cdots\times\CH(X_n)\times\CH(T)\times S.
$$
\end{lemma}

\begin{proof}
Since the restriction of $\CH(X)$ to the product is injective, the kernel of the composition
is the kernel of $\mathbb Z[x_{\alpha}]\rightarrow\CH(X)$. Since $\{x_{\alpha}\}$ is a generating
set for $\CH(X)$, we conclude that $\CH(X)=\mathbb Z[x_{\alpha}]/K$.
\end{proof}

We will also make frequent use of the excess intersection formula:

\begin{lemma}[Excess intersection]\label{excess intersection formula}
Consider a Cartesian diagram of smooth quotient stacks
$$
\begin{tikzcd}
Y' \arrow[r,"i'"] \arrow[d,"f'"'] & Y \arrow[d,"f"]\\
X' \arrow[r,"i"'] & X,
\end{tikzcd}
$$
where $f$ is proper and representable, and $i$ and $i'$ are regular local immersions.
Let the excess bundle $\mathcal E$ be
$$
0\rightarrow N_{Y'/Y}\rightarrow f'^*N_{X'/X}\rightarrow\mathcal E\rightarrow0.
$$
Then, for $z\in\CH(Y)$,
$$
i^*f_*z=f'_*\left(c_{\mathrm{top}}(\mathcal E)\cdot i'^*z\right).
$$
If $f$ and $f'$ are also regular local immersions, the same excess bundle
fits into the exact sequence
$$
0\rightarrow N_{Y'/X'}\rightarrow i'^*N_{Y/X}\rightarrow\mathcal E\rightarrow0.
$$
\end{lemma}

\begin{proof}
This follows from the compatibility of refined pullback with proper pushforward
\cite[Theorem 6.2(a)]{Fulton} and the excess intersection formula \cite[Theorem 6.3]{Fulton}.
\end{proof}

\subsection{Root gerbes and root stacks}
For a line bundle $L$ on a stack $X$, the \textit{root gerbe} $X(\sqrt[n]L)$ parametrizes
$n^{\text{th}}$ roots of $L$, and for an effective Cartier divisor $D\subseteq X$, the
\textit{root stack} $X(\sqrt[n]D)$ parametrizes $n^{\text{th}}$ roots of $\mathcal O_X(D)$ together
with an $n^{\text{th}}$ root of its canonical section.

\begin{lemma}\label{root gerbe Chow}
Let $X$ be a smooth quotient stack and $L$ a line bundle on $X$. Then
$$
\CH\left(X(\sqrt[n]L)\right)=\CH(X)[x]/(c_1(L)-nx),
$$
where $x$ is the first Chern class of the universal $n^{\text{th}}$ root of $L$.
\end{lemma}

\begin{proof}
Consider the complement of the zero section in the total space of $L$:
$U=\operatorname{Tot}(L)\setminus X$. Then $X(\sqrt[n]L)=[U/\mathbb G_m]$ where
$\mathbb G_m$ acts with weight $-n$. It then follows from excision that the Chow ring of
$X(\sqrt[n]L)$ is $\CH(X)[x]/(c_{\mathrm{top}}(L))$. Since $\mathbb G_m$ acts with weight $-n$, this
Chern class is $c_1(L)-nx$. Thus
$$
\CH\left(X(\sqrt[n]L)\right)=\CH(X)[x]/(c_1(L)-nx).
$$
\end{proof}

\begin{lemma}\label{root stack Chow}
Let $X$ be a smooth quotient stack and let $i:D\rightarrow X$ be a smooth effective Cartier divisor.
Suppose that $i^*:\CH(X)\rightarrow\CH(D)$ is surjective, and put $d=[D]$. Then
$$
\CH\left(X(\sqrt[n]D)\right)=\frac{\CH(X)[t]}{(nt-d, t\ker i^*)},
$$
where $t$ is the first Chern class of the universal root line bundle.
\end{lemma}

\begin{proof}
The root stack is the weighted blowup of $X$ along $D$ with weight $n$. Let $E$ be its exceptional
divisor. We write $s$ for the variable in \cite[Corollary 6.5]{AO}, so that $s=-[E]$.
Since $N_{D/X}=\mathcal O_X(D)|_D$, the equivariant top Chern class of the weighted normal
bundle is
$$
P(s)=i^*d+ns.
$$
The polynomial $Q$ in that corollary is therefore $Q(s)=ns+d$. Hence, using \cite[Theorem 7.1]{AO}
to extend to quotient stacks, we obtain
$$
\CH\left(X(\sqrt[n]D)\right)=\frac{\CH(X)[s]}{(ns+d, s\ker i^*)}.
$$
The universal root line bundle is $\mathcal O(E)$. Setting $t=-s=[E]=c_1(\mathcal O(E))$ proves
the lemma.
\end{proof}

\section{Plane quartics}\label{plane quartics section}

We open by considering the space $P$ of plane quartics and its locus of double conics, $Z$.
We deform $P$ to the normal cone of $Z$, and then find the Chow ring of the resulting stack
$X$ by using the patching lemma. After performing some excision in Sections
\ref{envelope section}--\ref{rigidification Chow ring}, this will serve as the base
for our root stack and root gerbe constructions of Section \ref{root section}.

\begin{definition}\label{blowup definition}
Let $V$ be a three-dimensional vector space, let
$P=\mathbb P(\Sym^4(V^\vee))$ be the space of plane quartics,
and let $Z=\{[q^2]:q\in\Sym^2(V^\vee)\}$ be the closed locus of double conics.
Let $Z_{\mathrm{sm}}$ and $Z_{\mathrm{sing}}$ be the loci of squares of smooth
and of singular conics, respectively. Set
$$
\mathcal D=\Bl_{Z\times0}(P\times\mathbb A^1_s)\setminus\widetilde{P\times0},\qquad
\mathcal D^\circ=\mathcal D\setminus\left(N_{Z/P}|_{Z_{\mathrm{sing}}}\right),
$$
where $\widetilde{P\times 0}$ is the strict transform. That is, $\mathcal D$ is the deformation of $P$
to the normal cone of $Z$, and $\mathcal D^\circ$ is obtained by deleting the part of its
special fiber that lies over singular conics. Let $\PGL_3\times\mathbb G_m$ act by
$$
([A],t)\cdot([F],s)=([F(A^{-1}X)],ts),\qquad
([A],t)\cdot([q^2],n)=([q(A^{-1}X)^2],t^{-1}[A]_*n),
$$
on the general and special fibers, respectively, where $A\in\GL_3$, $n\in N_{Z/P,[q^2]}$, and
$[A]_*$ denotes the induced action on the normal bundle $N_{Z/P}$. Write
$$
X=[\mathcal D^\circ/(\PGL_3\times\mathbb G_m)].
$$
Let $L$ be the line bundle associated with the character $(g,t)\mapsto t$ (i.e. the pullback of the
universal line bundle on $B\mathbb G_m$), and put $v=c_1(L)$.
The parameter $s$ is a section of $L$, and its zero locus is the special divisor $\widehat H$, so
$$
L=\mathcal O_{X}(\widehat H),\qquad N_{\widehat H/X}=L|_{\widehat H}.
$$
We call $[P/\PGL_3]$ the \textit{ordinary part} of $X$, and $\widehat H$ its \textit{special part},
or its \textit{special divisor}.
\end{definition}

\begin{definition}
If $q$ is a smooth conic in $\mathbb P(V)$, define
$B_{2k}=H^0(\mathbb P^1,\mathcal O_{\mathbb P^1}(2k))=H^0(q,\mathcal O_q(k))$, where,
for $q$, $\mathcal O_q(k)$ is being pulled back from $\mathbb P(V)$.
Note that $H^0(\mathbb P(V),\mathcal O_{\mathbb P(V)}(k))=\Sym^k(V^\vee)$.
On the special part, $\PGL_2$ acts on $V$ through its three-dimensional adjoint representation.
This identifies $\PGL_2$ with the stabilizer of $q$ in $\PGL_3$. The induced action on $B_{2k}$ is
$$
[A]\cdot F(X)=\det(A)^k F(A^{-1}X),\qquad A\in\GL_2.
$$
\end{definition}

\begin{lemma}\label{B_2k surjection lemma}
For $k\geq 2$, there is a surjection $\Sym^k(V^\vee)\rightarrow B_{2k}$ with kernel
$q\Sym^{k-2}(V^\vee)$.
\end{lemma}

\begin{proof}
The quadratic equation defining $q$ gives, for $k\geq2$, the exact sequence
$$
0\rightarrow\mathcal O_{\mathbb P(V)}(k-2)
\xrightarrow{\cdot q}\mathcal O_{\mathbb P(V)}(k)
\rightarrow\mathcal O_q(k)\rightarrow0.
$$
Since $H^1(\mathbb P(V),\mathcal O(k-2))=0$, taking global sections gives
$$
0\rightarrow\Sym^{k-2}(V^\vee)\xrightarrow{\cdot q}\Sym^k(V^\vee)\rightarrow B_{2k}\rightarrow0.
$$
\end{proof}

Before deleting the singular conics, the special fiber of $\mathcal D$ is the total space of
$N_{Z/P}$, on whose fibers $\mathbb G_m$ acts with weight $-1$. Therefore the special divisor is
the restriction of this vector bundle to the space of smooth conics. For a smooth conic $q$,
\cite[Proposition 4.1 and its proof]{Artebani} identifies the normal space $N_{Z/P,[q^2]}$ with $B_8$.
We will also rederive this in the proof of Lemma \ref{excess bundles}. Since $\PGL_3$ acts
transitively on smooth conics with stabilizer $\PGL_2$, we obtain
$$
\widehat H=\left[\frac{B_8}{\PGL_2\times\mathbb G_m}\right],
$$
where the action is given by $([A],t)\cdot F= t^{-1}([A]\cdot F)$. We write
$B_8^{\mathrm{sf}}\subseteq B_8$ for the open subset of sections with eight distinct zeros.

\begin{definition}
We write
$$
\mathcal F_4= \Sym^4(V^\vee)\otimes\det V,
$$
so that $A\in\GL_3$ acts by $(A\cdot f)(X)=\det(A) f(A^{-1}X)$.
\end{definition}

With this definition, a matrix of the form $tI$ acts as $t^{-1}$.
Therefore the central $\mathbb G_m$ of $\GL_3$ acts freely on $\mathcal F_4\setminus0$,
and the quotient by the $\mathbb G_m$ action is $P$, the space of plane quartics.
This induces the same action of $\GL_3/\mathbb G_m=\PGL_3$ as in Definition \ref{blowup definition}.
Taking the quotient therefore proves the following.

\begin{lemma}
We have the isomorphism
$$
[P/\PGL_3]\cong[(\mathcal F_4\setminus0)/\GL_3].
$$
\end{lemma}

\begin{definition}
Using the above identification, we will reuse $V$ to refer to the bundle associated with $V$ on
$[P/\PGL_3]$.
The projection $\mathcal D^\circ\rightarrow P$ induces a morphism $X\rightarrow[P/\PGL_3]$,
and we will also use $V$ to refer to the pullback of $V$ to $X$ along this map.
\end{definition}

\begin{note}\label{h_4 note}
When we work in $[P/\PGL_3]$, we will closely follow \cite{DLFV}. Using their same auxiliary
quotient by $\GL_3$, we have the identity
$$
\CH([P/\PGL_3])\cong \frac{\CH([P/\GL_3])}{(h_4-a)},
$$
where $a=c_1(V)$ and $h_4=c_1(\mathcal O_P(1))$ \cite[discussion before Theorem 7.2]{DLFV}.
Hence we will frequently make the substitution $h_4=a$.
\end{note}

Before we compute the Chow ring of the ordinary and special parts, we need one last fact.

\begin{lemma}[{\cite[Theorem 3]{Pandharipande}}]\label{PGL2 Chow}
We have
$$
\CH(B\PGL_2)=\mathbb Z[b,c]/(2c),
$$
where $b=c_2$ and $c=c_3$ are the Chern classes of the three-dimensional adjoint representation
of $\PGL_2$.
\end{lemma}

\begin{proof}
This follows using $\operatorname{SO}_3\cong\PGL_2$.
\end{proof}

\begin{lemma}\label{ordinary and special Chow rings}
The Chow ring of the ordinary part is
$$
\CH([P/\PGL_3])=\mathbb Z[a,b,c]/(P_{15}),
\qquad P_{15}=c_{15}(\Sym^4(V^\vee)\otimes\det V),
$$
where $(a,b,c)=(c_1,c_2,c_3)$ are the Chern classes of the standard representation
of $\GL_3$ and
$$
P_{15}=c(a^3-4ab+8c)(3a^3-16ab+64c)(36a^3c-8a^2b^2-182abc+36b^3+343c^2).
$$
The special divisor has Chow ring
$$
\CH(\widehat H)=\mathbb Z[b,c,v]/(2c),
$$
with $c_1(N_{\widehat H/X})=v$.
Lastly, the specialization map
$$
\CH([P/\PGL_3])\rightarrow\CH(\widehat H)/(v)
$$
appearing in the patching lemma is given by $a\mapsto0$, $b\mapsto b$, $c\mapsto c$.
In particular $P_{15}$ restricts to zero.
\end{lemma}

\begin{proof}
Since $\mathcal F_4$ is a vector bundle, $\CH([\mathcal F_4/\GL_3])=\mathbb Z[a,b,c]$.
Therefore, as $[P/\PGL_3]$ is obtained by excising the zero section of this rank 15 vector bundle,
$\CH([P/\PGL_3])=\mathbb Z[a,b,c]/(c_{15}(\mathcal F_4))$. Write $P_{15}$ for this Chern class,
and let $x_1,x_2,x_3$ be the Chern roots of $V$, so that
$$
a=x_1+x_2+x_3,\qquad
b=x_1x_2+x_1x_3+x_2x_3,\qquad
c=x_1x_2x_3.
$$
The roots of $\Sym^4(V^\vee)\otimes\det V$ are
$$
a-i_1x_1-i_2x_2-i_3x_3,\qquad i_1+i_2+i_3=4,
$$
and thus
$$
P_{15}=\prod_{\substack{i_1+i_2+i_3=4\\i_1,i_2,i_3\geq0}}(a-i_1x_1-i_2x_2-i_3x_3).
$$
\texttt{lemma-3-9.m2} computes this product.

For the special part, by the discussion after Lemma \ref{B_2k surjection lemma}, $\widehat H$ is a
vector bundle over $B\PGL_2\times B\mathbb G_m$ and hence has Chow ring
$\mathbb Z[b,c,v]/(2c)$ by Lemma \ref{PGL2 Chow}. Since, by construction,
$\widehat H$ is the vanishing of the weight-one parameter $s$ in
$$
\Bl_{Z\times 0}(P\times\mathbb A^1_s)\setminus\widetilde{P\times 0}
$$
we see that $c_1(N_{\widehat H/X})=v$.

The statement about $\CH([P/\PGL_3])\rightarrow\CH(\widehat H)/(v)$ follows because the standard
representation of $\GL_3$ restricts to the three-dimensional adjoint representation of $\PGL_2$.
Note that after setting $a=2c=0$, the factor $a^3-4ab+8c$ of $P_{15}$ vanishes and hence
so does $P_{15}$.
\end{proof}

\begin{proposition}\label{first patch}
There is an isomorphism
$$
\CH(X)=\frac{\mathbb Z[a,b,c,v]}{(P_{15}, av, 2vc)}.
$$
\end{proposition}

\begin{proof}
By Lemma \ref{ordinary and special Chow rings}, the top Chern class of the normal bundle of
$\widehat H$ is $v$, which isn't a zero-divisor in $\CH(\widehat H)$. Hence the patching lemma
(Lemma \ref{patching lemma}) applies and says that the Chow ring of $X$
is the fiber product of $\CH([P/\PGL_3])$ with $\CH(\widehat H)$ over $\CH(\widehat H)/(v)$.
Put
$$
R=\mathbb Z[a,b,c,v],\qquad I=(P_{15},v),\qquad J=(a,2c).
$$
Since $P_{15}\in J$ and $v$ isn't a zero-divisor in $R/J$, Lemma \ref{fiber product lemma} gives
$$
I\cap J=(P_{15})+vJ=(P_{15},av,2vc).
$$
Also $I+J=(a,2c,v)$, so the same lemma gives the claimed presentation.
\end{proof}

\section{A(n almost) Chow envelope for the boundary}\label{envelope section}
\begin{definition}
Let $J\subseteq X$ be the open substack consisting of smooth quartics (away from
$\widehat H$) and of nonzero squarefree binary octics (on $\widehat H$).
Let $H=J\cap\widehat H$ be its special divisor and $J_{\mathrm{nh}}=J\setminus H$.
\end{definition}

\begin{definition}
Let $B=X\setminus J$. Its ordinary part consists of all singular quartics. Its special part is
$$
\left[(B_8\setminus B_8^{\mathrm{sf}})/(\PGL_2\times\mathbb G_m)\right],
$$
the non-squarefree octics.
\end{definition}

We will define some important loci in $B$, and show that their Chow rings generate all cycles
pushed forward from $\CH(B)$ up to a factor of 3, \textit{almost} forming a Chow envelope.
We will then show in the following sections that in fact this factor of 3 poses no problem for us.

\subsection{The envelope}
\begin{definition}\label{D definitions}
For each of the following spaces, let $m_r:D_r\rightarrow P$ be the map which forgets the marking.
Each of the following definitions is from \cite{DLFV}.
\begin{enumerate}
\item Let $D_1\subseteq P\times\mathbb P(V)$ be the subscheme which parametrizes
pairs $(f,p)$ of singular quartics with a marked singular point. Thus $(f,p)\in D_1$ means
$f\in H^0(\mathbb P(V),I_p^2(4))$. See \cite[Remark 1.1 and Definition 3.4(1)]{DLFV}.
The projection $D_1\rightarrow\mathbb P(V)$ is a $\mathbb P^{11}$-bundle
by \cite[proof of Proposition 4.2]{FV18}.
\item Let $U=\Hilb^2(\mathbb P(V))$, let $\mathcal T\subseteq U\times\mathbb P(V)$
be its universal length-two subscheme, and let $\pi:U\times\mathbb P(V)\rightarrow U$ be
the projection. Let $\mathcal F=\pi_*I_{\mathcal T}^2(4)$. We define
$D_2=\mathbb P_U(\mathcal F)$. Its fiber over $T$ consists of quartics $f\in H^0(I_T^2(4))$.
See \cite[Proposition 5.2]{DLFV}.
\item Set $G=\Gr(2,\Sym^2(V^\vee))$ and write
$\mathcal S\subseteq\Sym^2(V^\vee)\otimes\mathcal O_G$
for its tautological bundle. Let
$D_{22}=\Hilb^2(\mathbb P(\Sym^2(V^\vee)))=\mathbb P_G(\Sym^2(\mathcal S^\vee))$
(see \cite[Lemma 4.7]{DLFV}). Its map to $P$ is Hilbert--Chow followed by multiplication of the
two conics.
\item $D_{13}=\mathbb P(\Sym^3(V^\vee))\times\mathbb P(V^\vee)$ maps by multiplication
of the cubic and the line.
\end{enumerate}
Let $W_r=D_r\times_PZ_{\mathrm{sm}}$, and put $\mathcal W_r=[W_r/\PGL_3]$. Deform $D_r$ to
the normal cone of its inverse image of $Z$, restrict over $\mathcal D^\circ$, and quotient by
$\PGL_3\times\mathbb G_m$. This defines $f_r:\mathcal Y_r\rightarrow X$. Note that the
ordinary part is $[D_r/\PGL_3]$.
\end{definition}

Each $f_r$ is proper and representable. Indeed, the maps $m_r$
are projective, and the induced map
$$
\Bl_{(D_r\times_P Z)\times0}(D_r\times\mathbb A^1_s)
\longrightarrow
\Bl_{Z\times0}(P\times\mathbb A^1_s)
$$
is proper. Hence, after quotienting by $\PGL_3\times\mathbb G_m$, the resulting maps $f_r$
are proper and representable.

Fix a smooth conic $q\subseteq\mathbb P(V)$. Since $\PGL_3$ acts transitively on the space
of smooth conics with stabilizer $\PGL_2$, it suffices to describe the fiber of each $W_r$ over $q^2$.

\begin{lemma}\label{excess bundles}
For $r=1,2,13,22$, the schemes $D_r$ and $W_r$ are smooth.
The fiber $F_r$ over a fixed doubled smooth conic $q^2$ is
$$
F_1=q,\qquad F_2=\Sym^2q\cong\mathbb P(V^\vee),\qquad F_{13}=\emptyset,\qquad
F_{22}=\mathbb P(B_4).
$$
Under the identification $N_{Z/P}|_{q^2}=B_8$,
we have the following excess sequences:
$$
\begin{gathered}
0\rightarrow B_6\otimes\mathcal O_q(-1)\rightarrow B_8\otimes\mathcal O_q\rightarrow\mathcal E_1
\rightarrow0,\\
0\rightarrow B_4\otimes\mathcal O_{\mathbb P(V^\vee)}(-2)\rightarrow B_8\otimes
\mathcal O_{\mathbb P(V^\vee)}\rightarrow\mathcal E_2\rightarrow0,\\
0\rightarrow\mathcal O_{\mathbb P(B_4)}(-2)\rightarrow B_8\otimes\mathcal O_{\mathbb P(B_4)}
\rightarrow\mathcal E_{22}\rightarrow0.
\end{gathered}
$$
The three inclusions are, respectively,
$$
s\otimes\ell_p\mapsto\ell_p s,\qquad s\otimes\ell^{\otimes2}\mapsto\ell^2s,\qquad
d^{\otimes2}\mapsto-d^2,
$$
where $\ell_p$ defines the tangent line to $q$ at $p$, $\ell$
defines the line cutting out $T\in\Sym^2q$, and $d\in B_4$
represents the conic direction.
\end{lemma}

\begin{proof}
The schemes $D_r$ are smooth by their projective bundle and product descriptions
in Definition \ref{D definitions}. We use the identification of the ambient normal bundle with $B_8$
from Section \ref{plane quartics section}.

The scheme $D_1$ parametrizes pairs of singular quartics and a marked singular point.
Since every point on $q^2$ is singular, the fiber of $D_1$ over $q^2$ is $F_1=q$ itself.
Similarly, since $D_2$ parametrizes pairs $(f,T)$ where $T$ has length 2 and
$f\in H^0(I_T^2(4))$, the fiber of $D_2$ over $q^2$ is $F_2=\Sym^2q$.
Intersecting with $q$ shows that $\Sym^2q\cong\mathbb P(V^\vee)$. Since $D_{13}$
parametrizes quartics which decompose as a cubic and a line, it has no fiber over $q^2$.
Lastly, consider $D_{22}$, which parametrizes pairs of conics, or, when the pair is the same
conic twice, a conic and a tangent direction. Its fiber over $q^2$ is thus the projectivization of the
tangent space at $q$: $F_{22}=\mathbb P(\Sym^2(V^\vee)/\langle q\rangle)=\mathbb P(B_4)$,
by Lemma \ref{B_2k surjection lemma}.

We now compute the normal bundles on these fibers. A first-order deformation of the quartic $q^2$
is $q^2+\epsilon F$, with $\epsilon^2=0$, where $F\in\Sym^4(V^\vee)$ is a plane quartic,
while deformations of it as a doubled conic have $F=2hq$,
$h\in\Sym^2(V^\vee)$, since $(q+\epsilon h)^2=q^2+2hq\epsilon$. Hence the normal bundle of
the special locus is $\Sym^4(V^\vee)/q\Sym^2(V^\vee)=B_8$, by Lemma \ref{B_2k surjection lemma}.

For $D_1$ and $D_2$, let $T\subseteq q$ denote the marking, of length
one or two, respectively, and put $t=(q^2,T)$. A tangent vector to
$D_r$ at $t$ consists of a deformation $q^2+\epsilon F$ together with
a compatible deformation $T_\epsilon$ of the marking. Restricting $F$
to $q$ gives a map
$$
\varphi: T_tD_r\rightarrow B_8,\qquad(F,T_\epsilon)\mapsto F|_q.
$$
We will compute its image and kernel.

First, we claim that the image is exactly
$$
H^0\left(q,\mathcal O_q(4)(-2T)\right)\subseteq B_8,
$$
that is, the global sections of $\mathcal O_q(4)$ which vanish on
$2T$. To see this, locally write $q=x$ and $I_T=(x,u)$, where $u$
cuts out $T$ on $q$. With this notation, a deformation of $T$ has
the form $I_{T_\epsilon}=(x+\epsilon a,u+\epsilon b)$. We can rewrite
the deformed quartic $x^2+\epsilon F$ as
$$
\begin{aligned}
x^2+\epsilon F&=\left((x+\epsilon a)-\epsilon a\right)^2+\epsilon F\\
&=(x+\epsilon a)^2+\epsilon(F-2ax).
\end{aligned}
$$
Hence the condition $x^2+\epsilon F\in I_{T_\epsilon}^2$ is equivalent to
$\epsilon(F-2ax)\in I_{T_\epsilon}^2$, which is equivalent to
$$
F-2ax\in I_T^2=(x^2,xu,u^2).
$$
Reducing modulo $x$, we see that $F|_q$ is divisible by $u^2$.
Thus every section in the image vanishes on $2T$.

Conversely, let $\beta\in H^0(q,\mathcal O_q(4)(-2T))$. Since the
restriction map $\Sym^4(V^\vee)\rightarrow B_8$ is surjective, we can lift $\beta$ to a
quartic $F$. Locally, its vanishing on $2T$ lets us write
$$
F=u^2c+xz.
$$
Choosing $a=z/2$ gives $F-2ax=u^2c\in I_T^2$. The choice of $a=z/2$ lifts
to a global deformation of $T$ because the sequence
$$
0\rightarrow N_{T/q}\rightarrow N_{T/\mathbb P(V)}
\rightarrow N_{q/\mathbb P(V)}|_T\rightarrow0
$$
is exact on global sections (since $T$ is finite). The resulting pair
$(q^2+\epsilon F,T_\epsilon)$ is a tangent vector to $D_r$ mapping
to $\beta$, which proves our claim about the image.

Next, we claim that the kernel is exactly $T_tW_r$. If $F|_q=0$, then $F=2qh$
for some quadratic $h$, and hence
$$
q^2+\epsilon F=(q+\epsilon h)^2.
$$
The compatibility condition above forces $a|_T=h|_T$, so the deformed marking lies on
$q+\epsilon h$. Thus the pair is in fact a deformation of the doubled conic,
and hence is tangent to $W_r$. Conversely, every deformation as a doubled conic has
$F|_q=0$. This proves our claim about the kernel.

Our descriptions of $F_1$ and $F_2$ give $\dim W_1=6$ and $\dim W_2=7$.
Since $\dim D_1=13$ and $\dim D_2=12$, the computations we just did
give $\dim T_tW_r=\dim W_r$ at every point, for $r=1,2$. Thus $W_1$ and $W_2$ are smooth.

Taking the quotient by $T_tW_r$, we conclude that
$$
(N_{W_r/D_r})_t\cong T_tD_r/T_tW_r\cong H^0\left(q,\mathcal O_q(4)(-2T)\right)\subseteq B_8.
$$

For $D_1$, write $T=p$. Let $\ell_p$ be an equation of the tangent line to $q$
at the marked point $p$. Its restriction to $q$ cuts out $2T$,
so the normal vectors are exactly $\ell_p g$, with $g\in B_6$.
The lines $\langle\ell_p\rangle$ form $\mathcal O_q(-1)$, giving
$$
\left.N_{W_1/D_1}\right|_{F_1}\cong B_6\otimes\mathcal O_q(-1),
\qquad g\otimes\ell_p\mapsto\ell_p g\in B_8.
$$

For $D_2$, write $T=q\cap\{\ell=0\}$. Vanishing on $2T$ is equivalent to divisibility by $\ell^2$,
so the normal vectors are exactly $\ell^2g$, with $g\in B_4$. Since $\ell$ spans
$\mathcal O(-1)$ on $F_2=\mathbb P(V^\vee)$, we obtain
$$
\left.N_{W_2/D_2}\right|_{F_2}\cong B_4\otimes\mathcal O_{\mathbb P(V^\vee)}(-2),
\qquad g\otimes\ell^{\otimes2}\mapsto\ell^2g\in B_8.
$$

Since $W_{22}\subseteq D_{22}$ is a divisor, its normal bundle is a line bundle. Hence at each
$[d]\in\mathbb P(B_4)$ it suffices to find a single vector in its span. The point $[d]$ specifies the
direction in which the doubled conic separates as $q+td$ and $q-td$. Their product is
$$
(q+td)(q-td)=q^2-t^2d^2=q^2-\epsilon d^2,\qquad \epsilon=t^2.
$$
Thus the normal bundle is $\langle d^2\rangle\subseteq B_8$.
Since $d$ spans $\mathcal O(-1)$ on $\mathbb P(B_4)$, this gives
$$
\left.N_{W_{22}/D_{22}}\right|_{F_{22}}
\cong\mathcal O_{\mathbb P(B_4)}(-2),
\qquad d^{\otimes2}\mapsto-d^2\in B_8.
$$
This also shows that $\dim T_{[d]}W_{22}=9=\dim W_{22}$, so $W_{22}$ is smooth.
\end{proof}

\begin{lemma}\label{special part isomorphisms}
Fix a smooth conic $q$, and let $B_0\subseteq\PGL_2$ be the
stabilizer of a point of $q$. Then
$$
\mathcal W_1= BB_0,\qquad
\mathcal W_2=[\mathbb P(V^\vee)/\PGL_2],\qquad
\mathcal W_{22}=[\mathbb P(B_4)/\PGL_2].
$$

For $r=1,2,22$, let $\mathcal N_r$ denote the vector bundle on $\mathcal W_r$ obtained by
descending the $\PGL_3$-equivariant normal bundle $N_{W_r/D_r}$.
Let $L$ be the universal line bundle on $B\mathbb G_m$. Then
$$
\mathcal Y_r^{\mathrm{sp}}=\operatorname{Tot}\left(\mathcal N_r\otimes L^{-1}\right)
\rightarrow\mathcal W_r\times B\mathbb G_m
$$
is a vector bundle. Consequently,
$$
\CH(\mathcal Y_r^{\mathrm{sp}})\cong\CH(\mathcal W_r)[v],\qquad v=c_1(L).
$$
\end{lemma}

\begin{proof}
Since $\PGL_3$ acts transitively on smooth conics with stabilizer $\PGL_2$, we have
$\mathcal W_r=[F_r/\PGL_2]$ where $F_r$ is the fiber over a smooth conic $q$, so the
description of each $\mathcal W_r$ follows from Lemma \ref{excess bundles}.

For $r=1,2,22$, the special fiber of the deformation to the normal cone over $W_r$ is
$N_{W_r/D_r}$. The action of $\mathbb G_m$ has weight $-1$, so taking the quotient yields
$$
\mathcal Y_r^{\mathrm{sp}}=\operatorname{Tot}\left(\mathcal N_r\otimes L^{-1}\right)
\rightarrow \mathcal W_r\times B\mathbb G_m.
$$
Vector bundle invariance therefore gives
$$
\CH(\mathcal Y_r^{\mathrm{sp}})=\CH(\mathcal W_r)[v].
$$
\end{proof}

We write $\eta=c_1(\mathcal O_{\mathbb P(B_4)}(1))$ on
$\mathcal W_{22}=[\mathbb P(B_4)/\PGL_2]$, and use the same notation
for its pullback to $\mathcal Y_{22}^{\mathrm{sp}}$.

\begin{lemma}\label{special restriction lemma}
For $r=1,2,22$, let
$$
i:\widehat H\rightarrow X,\qquad g:\mathcal W_r\times B\mathbb G_m\rightarrow B\PGL_2\times
B\mathbb G_m
$$
be the special-fiber inclusion and the morphism forgetting the
marking, respectively. For $z\in\CH(\mathcal Y_r)$, denote by
$z_{\mathrm{sp}}\in\CH(\mathcal Y_r^{\mathrm{sp}})=\CH(\mathcal W_r)[v]$ the special restriction
of $z$. Then
$$
i^*(f_r)_*z=g_*\left(z_{\mathrm{sp}} c_{\mathrm{top}}(\mathcal E_r\otimes L^{-1})\right)
\in\CH(\widehat H),
$$
where $\mathcal E_r$ is the excess bundle from Lemma \ref{excess bundles} and $L$ is the
universal line bundle on $B\mathbb G_m$.
\end{lemma}

\begin{proof}
Write $i_r:\mathcal Y_r^{\mathrm{sp}}\rightarrow\mathcal Y_r$ for the special-fiber inclusion and
$f_r^{\mathrm{sp}}$ for the restriction of $f_r$. Refined base change gives
$$
i^*(f_r)_*z=(f_r^{\mathrm{sp}})_*i_r^*z
$$
by \cite[Theorem 6.2(a) and (c)]{Fulton}.

By Lemma \ref{special part isomorphisms}, the special source and target are total spaces of their
normal bundles, twisted by $L^{-1}$. The square formed by $f_r^{\mathrm{sp}}$, $g$, and their zero
sections is Cartesian, and its excess bundle is $\mathcal E_r\otimes L^{-1}$ by Lemma
\ref{excess bundles}. The excess intersection formula (Lemma \ref{excess intersection formula})
then gives
$$
(f_r^{\mathrm{sp}})_*i_r^*z=g_*\left(z_{\mathrm{sp}} c_{\mathrm{top}}
(\mathcal E_r\otimes L^{-1})\right).
$$
\end{proof}

\subsection{Generators of $\CH(\mathcal Y_r)$}

Just as in Note \ref{h_4 note}, we have
$$
\CH([D_r/\PGL_3])=\CH([D_r/\GL_3])/(h_4-a),
$$
and we will frequently make the substitution $h_4=a$. We also write $V$ for the pullback of
$V$ to $\mathcal Y_r$.

\begin{lemma}\label{Y_1, Y_2, Y_22 generation}
The Chow rings of $\mathcal Y_1$, $\mathcal Y_2$, and $\mathcal Y_{22}$
are generated, as modules over $\CH(X)$, by the following elements:
\begin{enumerate}
\item on $\mathcal Y_1$, the classes $1,\rho,\rho^2$, where
$\rho=c_1(\mathcal O_{\mathbb P(V)}(1))$;
\item on $\mathcal Y_2$, the classes $h_1^is^j$, $0\leq i,j\leq2$,
where $h_1$ is the pullback of $c_1(\mathcal O_{\mathbb P(V^\vee)}(1))$ under
$$
U\rightarrow\mathbb P(V^\vee),\qquad T\mapsto\ell_T,
$$
with $\ell_T$ the unique line containing $T$, and
$s=c_1(\mathcal L^\vee)$, where $\mathcal L$ is the line bundle
on $U$ with fiber
$$
\mathcal L_T=H^0(\ell_T,\mathcal I_{T/\ell_T}(2));
$$
\item on $\mathcal Y_{22}$, all monomials in $s_1=c_1(\mathcal S)$,
$s_2=c_2(\mathcal S)$, and
$\tau=c_1(\mathcal O_{\mathbb P_G(\Sym^2(\mathcal S^\vee))}(1))$.
\end{enumerate}
\end{lemma}

\begin{proof}
We first analyze the ordinary parts.
For $D_1$, apply \cite[Proposition 4.1 and the proof of Proposition 4.2]{FV18} with $(n,d)=(3,4)$.
Their $\widetilde Z$ is our $D_1$, and their classes $(h_d,t)$ are $(h_4,\rho)$.
This shows that $1,\rho,\rho^2$ generate over $\CH(X)$. For $D_2$,
\cite[Lemma 5.1 and Proposition 5.2]{DLFV} gives the generators $h_1^is^j$, $0\leq i,j\leq2$.
For $D_{22}$, use \cite[Lemma 4.7 and the first paragraph of the
proof of Proposition 4.8]{DLFV}. Their generators involve $\tau$ and the Chern classes of the
tautological quotient on $G$. The identity
$$
c_t\left( \Sym^2(V^\vee)/\mathcal S\right)=\frac{c_t( \Sym^2(V^\vee))}{1+s_1t+s_2t^2}
$$
expresses those Chern classes as polynomials in $a,b,c,s_1,s_2$.
Therefore the claimed generators suffice for the ordinary part.

Let $i_r:\mathcal Y_r^{\mathrm{sp}}\rightarrow\mathcal Y_r$ be the inclusion.
By Lemma \ref{special part isomorphisms},
$$
\CH(\mathcal Y_r^{\mathrm{sp}})=\CH(\mathcal W_r)[v],
$$
and we will show that $\CH(\mathcal W_r)$ is generated by the claimed elements
as a module over $\CH(B\PGL_2)=\mathbb Z[b,c]/(2c)$.

For $r=1$, Lemma \ref{special part isomorphisms} says $\mathcal W_1= BB_0$ where
the Borel subgroup $B_0\subseteq\PGL_2$ is the stabilizer of a point $p$ on a smooth conic $q$.
Choose coordinates in which $p=[1:0:0]$ and a maximal torus $T\subseteq B_0$ acts on
$V$ by $\operatorname{diag}(t^{-1},1,t)$.
The restriction of $\mathcal O_{\mathbb P(V)}(-1)$ at $p$ has character $t^{-1}$,
so the restriction of $\rho=c_1(\mathcal O_{\mathbb P(V)}(1))$ is the first Chern class of the
generating character $t$. By \cite[Lemma 3]{EG97}, restriction to $T$ therefore gives
$$
\CH(\mathcal W_1)=\CH(BT)=\mathbb Z[\rho].
$$
Since the characters are $t^{-1},1,t$, the Chern roots of $V|_{BT}$ are
$-\rho,0,\rho$, so this restriction isomorphism also gives
$b=c_2(V)=-\rho^2$ on $\mathcal W_1$. Hence $\CH(\mathcal W_1)$ is generated as a
$\CH(B\PGL_2)$-module by $1,\rho,\rho^2$ (in fact by $1,\rho$).

For $r=2$, Lemma \ref{special part isomorphisms} gives $\mathcal W_2=[\mathbb P(V^\vee)/\PGL_2]$,
and so the projective bundle formula says that the Chow ring is generated by $1,h_1,h_1^2$
over $\CH(B\PGL_2)$.
We will also determine the special restriction of $s$, since it will be useful later. For
$T\in U=\Hilb^2(\mathbb P(V))$, let $\ell$ be the line it determines. The fiber of $\mathcal L$ at
$T$ is
$$
H^0(\ell,\mathcal I_{T/\ell}(2)).
$$
On $F_2$, we have $T=q\cap\ell$, so this line is spanned by $q|_\ell$.
Since $q$ contains no line, these restrictions are everywhere nonzero, and hence give a
nonvanishing section. That is, this line bundle is trivial, which shows
\begin{equation}\label{s vanishing}
s=0\quad\text{on }\mathcal W_2.
\end{equation}

For $r=22$, Lemma \ref{special part isomorphisms} gives
$$
\mathcal W_{22}=[\mathbb P(B_4)/\PGL_2].
$$
Recalling that  $\eta=c_1(\mathcal O_{\mathbb P(B_4)}(1))$, the projective bundle formula says
that the Chow ring is generated by $1,\eta,\dots,\eta^4$ over $\CH(B\PGL_2)$, and so all we have
to do is express $\eta$ in terms of our desired generators. Under the identification
$\mathbb P(B_4)=\mathbb P(\Sym^2(V^\vee)/\langle q\rangle)$ from Lemma
\ref{B_2k surjection lemma}, projection modulo $q$ gives
$$
\mathcal S|_{\mathcal W_{22}}/(\mathcal O\cdot q)=\mathcal O_{\mathbb P(B_4)}(-1).
$$
Thus, on $\mathcal W_{22}$,
\begin{equation}\label{s_2 vanishing}
s_1=c_1(\mathcal S)=c_1(\mathcal O)+c_1(\mathcal O(-1))=-\eta,\qquad
s_2=c_2(\mathcal S)=0.
\end{equation}
This shows that the Chow ring is generated by the claimed generators over $\CH(B\PGL_2)$, but
we will also determine the special restriction of $\tau$, since this will be useful later.
Dualizing gives
$$
\Ann(q)=\{\varphi\in\mathcal S^\vee:\varphi(q)=0\}\cong(\mathcal S/\langle q\rangle)^\vee=
\mathcal O_{\mathbb P(B_4)}(1).
$$
A point of $\mathcal W_{22}$ corresponds to the double point at $[q]$ in $\mathbb P(\mathcal S)$.
If $\varphi$ spans $\Ann(q)$, this double point is cut out by $\varphi^2$. Since
$D_{22}=\mathbb P_G(\Sym^2(\mathcal S^\vee))$, the fiber of its tautological line bundle is precisely
the line spanned by $\varphi^2$. Thus, on $\mathcal W_{22}$, we have
$$
\mathcal O_{D_{22}}(-1)=\Ann(q)^{\otimes 2}=\mathcal O_{\mathbb P(B_4)}(2).
$$
Taking the dual says that the special restriction of $\tau$ is
\begin{equation}\label{tau formula}
\tau=-2\eta.
\end{equation}

For each $r$, let $M_r\subseteq\CH(\mathcal Y_r)$ be the $\CH(X)$-span of our claimed generators.
The preceding calculations show that $M_r$ surjects onto the Chow rings of both the ordinary and
special parts of $\mathcal Y_r$. The special divisor is the pullback of $\widehat H$, so by
Definition \ref{blowup definition},
$$
(i_r)_*1=f_r^*v.
$$
Lemma \ref{Chow generation lemma} therefore gives $M_r=\CH(\mathcal Y_r)$.
\end{proof}

\begin{theorem}\label{degree three}
The maps $f_r:\mathcal Y_r\rightarrow B$ satisfy
$$
3\CH(B)\subseteq\sum_{r=1,2,22}\im\left((f_r)_*:\CH(\mathcal Y_r)\rightarrow\CH(B)\right).
$$
\end{theorem}

\begin{proof}
By using an equivariant approximation, we can assume that all relevant spaces are schemes.
Let $Y\subseteq B$ be an integral subvariety with function field $K$. We will find an integral subvariety
$\widetilde Y\subseteq\mathcal Y_r$, for some $r$, such that $\widetilde Y\rightarrow Y$ has degree
one or three. Let $G=\operatorname{Gal}(\overline K/K)$.

Suppose first that the generic point of $Y$ lies in the ordinary part, so that it represents a singular
quartic in $\mathbb P(V)$. This quartic is either geometrically integral, geometrically reduced but
reducible, or geometrically nonreduced.

Assume it is geometrically integral. Then it has at most three singular points over $\overline K$.
If it has only one, then the action of $G$ fixes it and hence the point is defined over $K$.
This gives a $K$-point of $\mathcal Y_1$.
If it has two points over $\overline K$, then the unordered pair of these points is fixed by $G$
and hence defined over $K$. This defines a length-two scheme over $K$,
and hence a $K$-point of $\mathcal Y_2$. If there are three singular points,
then the action of $G$ either fixes a point or acts transitively on all three. If a point is fixed,
then this gives a $K$-point of $\mathcal Y_1$. If it acts transitively, then the collection of all
three points is fixed by $G$ and hence defines an $L$-point of $\mathcal Y_1$ for
$[L:K]=3$.

Now assume it is geometrically reduced but not irreducible.
If it is a (geometrically irreducible) cubic and line, then each irreducible component is fixed by $G$
and hence is necessarily defined over $K$. The intersection is then length three, which necessarily
gives at least one $K$-point of $\mathcal Y_1$ or an $L$-point for $[L:K]=3$. If it is two conics
(one of them possibly non-irreducible),
then the pair is fixed by $G$ and also must be defined over $K$,
giving a $K$-point of $\mathcal Y_{22}$. If it is four distinct lines, consider the collection of the
three unordered pairings of lines. Either one of these pairs is fixed by the Galois action,
in which case it gives a $K$-point of $\mathcal Y_{22}$, or the action of $G$ is transitive on the
collection, in which case it defines an $L$-point for $[L:K]=3$.

Now assume it is nonreduced. If a line is repeated and defined over $K$, then every point on it is
singular. Picking any of these defines a $K$-point of $\mathcal Y_1$. If a line is repeated and
not defined over $K$, let $\ell$ be the line. Since $\ell$ is not defined over $K$, its
orbit under $G$ must contain exactly one other line. Intersecting that line with $\ell$
gives a $G$-invariant singular $K$-point and hence a $K$-point of $\mathcal Y_1$. Now assume
it is a doubled smooth conic, $q^2$. Then $q$ is necessarily defined over $K$ since it is fixed by $G$.
Intersecting $q$ with a line gives a length-two scheme over $K$ and hence a $K$-point
of $\mathcal Y_2$.

This covers all cases where the generic point lies in the ordinary part, so now we assume
that it lies in the special part. That is, we have a smooth conic $q$ over $K$ and a
non-squarefree section $F$ of $\mathcal O_q(4)$. If $F\neq 0$, define
$$
T=\sum_{p\in q(\overline K)}\left\lfloor\frac{\operatorname{mult}_p\ddiv(F)}2\right\rfloor p,
\qquad m=\deg T.
$$
Note that $1\leq m\leq 4$. Since $T$ is fixed by $G$, it is defined over $K$.
Arguing as above, $m=1$ gives a $K$-point of $\mathcal Y_1$,
$m=2$ gives a $K$-point in $\mathcal Y_2$, and $m=3$ gives a point possibly of degree one or three
in $\mathcal Y_1$. If $m=4$, then $\ddiv(F)=2T$.
Writing $T=\ddiv(d)$ with $d\in H^0(q,\mathcal O_q(2))$, we have $F=\lambda d^2$,
which gives a $K$-point in $\mathcal Y_{22}$. If $F=0$, intersect $q$ with any line
to get a $K$-point of $\mathcal Y_2$.

In each case, we have found a point over the generic point of $Y$ with residue field of degree one
or three over $K$. Taking its closure in $\mathcal Y_r$ gives the desired $\widetilde Y$, and
$(f_r)_*[\widetilde Y]$ is either $[Y]$ or $3[Y]$. Thus $3[Y]$ lies in the image of pushforward,
as desired.
\end{proof}

\section{The boundary relations}\label{pushforward section}
We now push forward the generators of each $\mathcal Y_r$, constructing a list
of relations on $\CH(X)$ which, in the next section, we will show is complete. In the definition below,
the polynomials $A_1,A_2,A_3$ will later be shown to be relations pushed forward from
$\mathcal Y_1$, $A_4$ will come from $\mathcal Y_{13}$, and $A_5,A_6$ will come from
$\mathcal Y_2$.

\begin{definition}\label{candidate ring}
In $\mathbb Z[a,b,c,v]$ put
$$
\begin{aligned}
A_1&=-9a-14v,& A_2&=12a^2-24b+2v^2,\\
A_3&=-12a^3+34ab-28c+14bv,& A_4&=19a^3-88ab+224c,\\
A_5&=3a^3-26ab+308c+50vb-10v^3,\\
A_6&=3a^4-23a^2b+86ac+24b^2-vc-35v^2b+v^4.
\end{aligned}
$$
Set
$$
I_{\mathrm{exc}}=(P_{15},av,2vc,A_1,A_2,A_3,A_4,A_5,A_6)
\subseteq\mathbb Z[a,b,c,v],
$$
and put
$$
Q=\mathbb Z[a,b,c,v]/I_{\mathrm{exc}}.
$$
Note that there is a natural map $\CH(X)\rightarrow Q$, since $P_{15},av,2vc$ are exactly
the relations on $\CH(X)$ by Proposition \ref{first patch}.
\end{definition}

\begin{lemma}\label{112 annihilator}
Restriction to the ordinary smooth quartics gives
$$
\CH(J_{\mathrm{nh}})=Q/(v).
$$
Moreover, $112 \omega=0$ for every $\omega\in(v)\subseteq Q$ of degree at least two.
\end{lemma}

\begin{proof}
Theorem 7.1 of \cite{DLFV} shows that the locus of smooth quartics in
$[P/\GL_3]$ has the relations $\alpha_1,\alpha_2,\alpha_3,\delta_{\{1,3\}}$.
Making the substitutions $c_1=a,c_2=b,c_3=c,h_4=a$ gives
$-9a$, $3a^2-24b$, $ab-28c$, and $A_4$. Thus
$$
\CH(J_{\mathrm{nh}})=\mathbb Z[a,b,c]/(9a,3a^2-24b,ab-28c,A_4).
$$
The script \texttt{lemma-5-2.m2} checks that $P_{15}$, $A_5|_{v=0}$, and
$A_6|_{v=0}$ belong to this ideal, which proves the first claim.
For the second, it checks
$$
112av,\quad112bv,\quad112cv,\quad112v^2\in I_{\mathrm{exc}}.
$$
\end{proof}

\begin{lemma}\label{excess chern classes}
For the excess bundles of Lemma \ref{excess bundles}, we have
$$
\begin{aligned}
c_2(\mathcal E_1\otimes L^{-1})
&=(3\rho-v)(4\rho-v),\\[2pt]
c_4(\mathcal E_2\otimes L^{-1})
&=v^4-10h_1v^3+(60h_1^2+25b)v^2\\
&\qquad\qquad +(50bh_1+c)v+120bh_1^2+144b^2,\\[2pt]
c_8(\mathcal E_{22}\otimes L^{-1})
&=(v^2+b)(v^2+9b)\bigl[
(v^2+4b)(v^2+16b)\\
&\qquad\qquad
+2\eta(2\eta-v)(v^2+4\eta^2+20b)
\bigr]+c^2v^2.
\end{aligned}
$$
\end{lemma}

\begin{proof}
Lemma \ref{B_2k surjection lemma} gives
$$
0\rightarrow\Sym^{k-2}(V^\vee)\rightarrow\Sym^k(V^\vee)\rightarrow B_{2k}\rightarrow0.
$$
Consequently,
$$
c_t(B_{2k})=\frac{c_t(\Sym^k(V^\vee))}{c_t(\Sym^{k-2}(V^\vee))}.
$$
For formal Chern roots $x,y,z$ of $V$, we have
$$
c_t(\Sym^m(V^\vee))=\prod_{i+j+\ell=m}\left(1-(ix+jy+\ell z)t\right).
$$

For a rank-$r$ bundle $E$ with formal Chern roots $l_1,\dots,l_r$, and a line bundle
$M$ with $u=c_1(M)$,
$$
\begin{aligned}
c_t(E\otimes M)&=\prod_{j=1}^r\left(1+(l_j+u)t\right)\\
&=\sum_{i=0}^r c_i(E)t^i(1+ut)^{r-i}.
\end{aligned}
$$
So, putting
$$
\Phi_{2k}(t,u)=
\sum_{i=0}^{2k+1}c_i(B_{2k})t^i(1+ut)^{2k+1-i},
$$
we see that $\Phi_{2k}(t,u)=c_t(B_{2k}\otimes M)$ where $M$ is any line bundle with
$c_1(M)=u$. Our above expression for $c_t(B_{2k})$ in terms of $c_t(\Sym^k(V^\vee))$
allows us to explicitly compute $\Phi_{2k}(t,u)$.

Twisting the three sequences of Lemma \ref{excess bundles}
by $L^{-1}$, their middle terms are $B_8\otimes L^{-1}$ and their left-hand terms are, respectively,
$$
B_6\otimes\mathcal O(-1)\otimes L^{-1},\qquad
B_4\otimes\mathcal O(-2)\otimes L^{-1},\qquad
\mathcal O(-2)\otimes L^{-1}.
$$
The twisting line bundles have first Chern classes
$-\rho-v$, $-2h_1-v$, and $-2\eta-v$, respectively.
Therefore
$$
\begin{aligned}
c_t(\mathcal E_1\otimes L^{-1}) & =\frac{\Phi_8(t,-v)}{\Phi_6(t,-\rho-v)},\\
c_t(\mathcal E_2\otimes L^{-1}) & =\frac{\Phi_8(t,-v)}{\Phi_4(t,-2h_1-v)},\\
c_t(\mathcal E_{22}\otimes L^{-1}) & =\frac{\Phi_8(t,-v)}{1+(-2\eta-v)t}.
\end{aligned}
$$
Carrying out this calculation in \texttt{lemma-5-3.m2} and extracting the degree $2,4,8$
terms proves the lemma.
\end{proof}

Our process for determining the pushforwards of elements from $\mathcal Y_r$ to $X$ will
proceed more or less the same in each of the following propositions. We will determine the
ordinary and special restrictions of the pushforwards, and then use the following lemma
to determine the pushforward using these restrictions.

\begin{lemma}\label{reconstruction lemma}
Let $z\in\CH(X)$ have ordinary and special restrictions represented by
$O(a,b,c)\in\mathbb Z[a,b,c]$ and $S(b,c,v)\in\mathbb Z[b,c,v]$, respectively.
Then
$$
z=O(a,b,c)+S(b,c,v)-S(b,c,0)\quad\text{in }\CH(X).
$$
\end{lemma}

\begin{proof}
By the fiber-product description in the proof of Proposition \ref{first patch},
the two restrictions determine $z$. By construction of fiber products, their images in
$\CH(\widehat H)/(v)=\mathbb Z[b,c]/(2c)$ must agree. That is, we must have
$$
O(0,b,c)=S(b,c,0)\quad\text{in }\mathbb Z[b,c]/(2c).
$$
Hence the claimed expression for $z$ restricts to $O(a,b,c)$ on the ordinary part,
where $v=0$, and to
$$
O(0,b,c)+S(b,c,v)-S(b,c,0)=S(b,c,v)
$$
on the special part, where $a=0$. It therefore equals $z$.
\end{proof}

\begin{proposition}\label{Y_1 and Y_13 pushforwards}
The pushforwards of $1,\rho,\rho^2$ along $f_1:\mathcal Y_1\rightarrow X$ are
$$
A_1,\qquad A_2,\qquad A_3.
$$
The fundamental class of $\mathcal Y_{13}$ pushes forward to $A_4$.
\end{proposition}

\begin{proof}
For the ordinary restriction, let
$$
\pi:P\times\mathbb P(V)\rightarrow P
$$
be the projection. By \cite[Proposition 4.3 and Lemma 4.4]{FV18}, with $n=3$ and $d=4$,
the ordinary restriction of the pushforward of $\rho^k$ is obtained by reducing
$$
\rho^k\left((h_4+3\rho)^3-a(h_4+3\rho)^2+b(h_4+3\rho)-c\right)
$$
using the projective bundle relation
$$
\rho^3+a \rho^2+b \rho+c=0,
$$
then extracting the coefficient of $\rho^2$ and making the substitution $h_4=a$ from
Note \ref{h_4 note}. Here their $c_1,c_2,c_3,t$ are our $a,b,c,\rho$. Carrying this out in
\texttt{proposition-5-5.m2} gives
$$
\begin{aligned}
1&\mapsto -9a,\\
\rho&\mapsto 12a^2-24b,\\
\rho^2&\mapsto -12a^3+34ab-28c.
\end{aligned}
$$

By Lemmas \ref{special restriction lemma} and \ref{excess chern classes}, we know that
the special restriction of the pushforward of $\rho^k$ is
$$
g_*\left(\rho^k c_2(\mathcal E_1\otimes L^{-1})\right)=g_*\left(\rho^k(4\rho-v)(3\rho-v)\right),
$$
where
$$
g:\mathcal W_1\times B\mathbb G_m=[q/\PGL_2]\times B\mathbb G_m\rightarrow
B\PGL_2\times B\mathbb G_m.
$$
By the proof of Lemma \ref{Y_1, Y_2, Y_22 generation},
we have $\rho^2=-g^*b$.
Since $g$ has relative dimension one and
$\mathcal O_q(1)$ has degree two,
$$
g_*1=0,\qquad g_*\rho=2.
$$
The projection formula therefore gives
$$
g_*(\rho^{2m})=0,\qquad g_*(\rho^{2m+1})=2(-b)^m\qquad(m\geq0),
$$
which yields
$$
1\mapsto-14v,\qquad
\rho\mapsto-24b+2v^2,\qquad
\rho^2\mapsto14bv.
$$

Combining the ordinary and special restrictions using Lemma \ref{reconstruction lemma}
proves the proposition for $\mathcal Y_1$: we obtain
$$
(f_1)_*1=A_1,\qquad (f_1)_*\rho=A_2,\qquad (f_1)_*\rho^2=A_3.
$$

Finally, \cite[\S 4.3, before Lemma 4.9]{DLFV} gives the ordinary restriction of the fundamental
class of $\mathcal Y_{13}$. Using our notation, it is $A_4=19a^3-88ab+224c$.
Since $\mathcal Y_{13}$ doesn't have a special fiber, we conclude that its fundamental class in
$\CH(X)$ is $A_4$.
\end{proof}

We now compute the pushforwards of the generators of $\CH(\mathcal Y_2)$.

\begin{proposition}\label{Y_2 pushforward}
The pushforwards of $1,h_1,h_1^2$ from $\mathcal Y_2$ are
$$
-5(a+v)A_1-5A_2,\qquad A_5,\qquad A_6.
$$
The pushforwards of $h_1^is,h_1^is^2$ contribute no new relations.
\end{proposition}

\begin{proof}
We first consider $1,h_1,h_1^2$.

On the ordinary part, $D_2$ plays the role of $D$ from \cite[Proposition 5.2]{DLFV},
the ordinary restriction of the pushforward of $h_1^i$ is their $p_*q^*h_1^is^0$, and, since
$l_1,l_2,l_3$ are the Chern roots of $V^\vee$ \cite[\S 2]{DLFV},
$$
h_4=a,\qquad a=-(l_1+l_2+l_3),\qquad b=l_1l_2+l_1l_3+l_2l_3,\qquad c=-l_1l_2l_3.
$$
Evaluating the formula in \cite[end of the proof of Proposition 5.4]{DLFV} with these
substitutions, \texttt{proposition-5-6.m2} gives
$$
\begin{aligned}
1&\mapsto-15a^2+120b,\\
h_1&\mapsto3a^3-26ab+308c,\\
h_1^2&\mapsto3a^4-23a^2b+86ac+24b^2.
\end{aligned}
$$

Since $\mathcal W_2=\mathbb P(V^\vee)$ over $B\PGL_2$, we have
$$
g_*1=0,\qquad g_*h_1=0,\qquad g_*h_1^2=1
$$
and
$$
h_1^3+bh_1-c=0.
$$
From Lemma \ref{excess chern classes}, we have an explicit expression for
$c_4(\mathcal E_2\otimes L^{-1})$. We plug this into Lemma \ref{special restriction lemma}
and use our above pushforwards of powers of $h_1$ to compute
$$
g_*\left(h_1^ic_4(\mathcal E_2\otimes L^{-1})\right).
$$
This computation is carried out in \texttt{proposition-5-6.m2}, obtaining, respectively,
$$
120b+60v^2,\qquad 50vb-10v^3,\qquad24b^2-vc-35v^2b+v^4.
$$
Combining this with the ordinary restriction using Lemma \ref{reconstruction lemma} proves the
proposition for $1,h_1,h_1^2$.

For the ordinary part of $h_1^is,h_1^is^2$, we again use the formula from
\cite[end of the proof of Proposition 5.4]{DLFV}. Since the special restriction of $s$ is 0
(see \eqref{s vanishing}), the special restrictions of all of these elements vanish.
Hence the pushforwards are just the ordinary restrictions, and \texttt{proposition-5-6.m2} verifies that
these restrictions don't add any new relations.
\end{proof}

Finally, we compute the pushforwards of the generators of $\CH(\mathcal Y_{22})$, and see
that they contribute no new relations.

\begin{proposition}\label{no relations}
The composition
$$
\CH(\mathcal Y_{22})\xrightarrow{(f_{22})_*}\CH(X)\rightarrow Q
$$
is zero.
\end{proposition}

\begin{proof}
Using projective bundle and Grassmannian bundle relations, we only need to consider
$s_1^ps_2^q\tau^r$ for $p,q\leq 4$ and $r\leq 2$.

For the ordinary restrictions, we use the localization computation in
\cite[proof of Proposition 4.8]{DLFV}. Their map $f$ is our map
$m_{22}$. Let $O_{p,q,r}$ be the ordinary restriction of the pushforward of
$s_1^ps_2^q\tau^r$. We use their same notation for $T,Q_i,\chi_i$ for the diagonal
$\mathbb G_m^3\subseteq\GL_3$, the six quadratic monomials, and their characters.

Applying their localization formula to $s_1^ps_2^q\tau^r$ gives
$$
\begin{aligned}
O_{p,q,r}&=\sum_{(Q_i,Q_j)}\frac{(s_1^ps_2^q\tau^r)|_{(Q_i,Q_j)}}
{c_{\mathrm{top}}^T(T_{(Q_i,Q_j)}D_{22})}[m_{22}((Q_i,Q_j))]_T\\
&\quad+\sum_{(Q_i>Q_j)}\frac{(s_1^ps_2^q\tau^r)|_{(Q_i>Q_j)}}
{c_{\mathrm{top}}^T(T_{(Q_i>Q_j)}D_{22})}[m_{22}((Q_i>Q_j))]_T.
\end{aligned}
$$
The first sum runs over the fifteen reduced pairs $i<j$, and the second over the thirty
nonreduced $i\neq j$ (supported at $Q_i$ in the direction of $Q_j$). We reuse their explicit
computations of each of these terms, except for $\tau|_{(Q_i>Q_j)}$, which contains a small typo:
it is equal to $2\chi_j$, since, again using their notation, the subscheme of $F_{ij}$ of length two
supported at $Q_i$ is the vanishing of the form $a_j^2$. This then gives
$$
\begin{aligned}
\left.s_1\right|_{(Q_i,Q_j)}&=\left.s_1\right|_{(Q_i>Q_j)}=\chi_i+\chi_j,\\
\left.s_2\right|_{(Q_i,Q_j)}&=\left.s_2\right|_{(Q_i>Q_j)}=\chi_i\chi_j,\\
\left.\tau\right|_{(Q_i,Q_j)}&=\chi_i+\chi_j,\qquad \left.\tau\right|_{(Q_i>Q_j)}=2\chi_j,\\
c_{\mathrm{top}}^T(T_{(Q_i,Q_j)}D_{22})&=\prod_{k\neq i}(\chi_k-\chi_i)
\prod_{k\neq j}(\chi_k-\chi_j),\\
c_{\mathrm{top}}^T(T_{(Q_i>Q_j)}D_{22})&=2(\chi_i-\chi_j)^2
\prod_{k\neq i,j}(\chi_k-\chi_i)(\chi_k-\chi_j),\\
m_{22}((Q_i,Q_j))&=[Q_iQ_j],\qquad m_{22}((Q_i>Q_j))=[Q_i^2],\\
[m_{22}((Q_i,Q_j))]_T&=\prod_{F\neq Q_iQ_j}(a+\chi_F),\\
[m_{22}((Q_i>Q_j))]_T&=\prod_{F\neq Q_i^2}(a+\chi_F),
\end{aligned}
$$
where $\chi_F$ is the character of $F$ (their $\theta_i$ when $F=F_i$).
\texttt{proposition-5-7.m2} evaluates the formula for $O_{p,q,r}$ using
$$
\begin{aligned}
\chi_1&=2l_1, & \chi_2&=l_1+l_2, & \chi_3&=l_1+l_3,\\
\chi_4&=2l_2, & \chi_5&=l_2+l_3, & \chi_6&=2l_3,
\end{aligned}
$$
where the $l_i$ are the Chern roots of $V^\vee$, that is, where
$$
a=-(l_1+l_2+l_3),\qquad b=l_1l_2+l_1l_3+l_2l_3,\qquad c=-l_1l_2l_3.
$$

On the special fiber, we have $s_1=-\eta$, $s_2=0$, $\tau=-2\eta$ (see \eqref{s_2 vanishing}
and \eqref{tau formula}).
Thus the special restriction of $s_1^ps_2^q\tau^r$ is zero when $q>0$.
In the $q=0$ case, Lemma \ref{special restriction lemma} says that the special restriction of the
pushforward is
$$
S_{p,r}:=g_*\left((-\eta)^p(-2\eta)^rc_8(\mathcal E_{22}\otimes L^{-1})\right),
\qquad g_*\eta^{4+j}=s_j(B_4),
$$
where $s_j(B_4)$ is the Segre class. Set $S_{p,q,r}=S_{p,r}$ if $q=0$ and zero otherwise.
In \texttt{proposition-5-7.m2} we use the explicit expression for $c_8(\mathcal E_{22}\otimes L^{-1})$
and the Chern classes of $B_4$ to carry out the computation of $S_{p,q,r}$.

Combining the ordinary and special parts using Lemma \ref{reconstruction lemma},
\texttt{proposition-5-7.m2} shows that each pushforward vanishes in $Q$.
\end{proof}

\section{Completing the excision}\label{rigidification Chow ring}
We now complete the excision of $B$ and collect some results which will be of use later.
Recall our ring $Q$ from Definition \ref{candidate ring}:
$$
Q=\frac{\mathbb Z[a,b,c,v]}{(P_{15},av,2vc,A_1,A_2,A_3,A_4,A_5,A_6)}.
$$

\begin{proposition}\label{excision result}
We have
$$
\CH(J)=Q.
$$
\end{proposition}

\begin{proof}
Recall that $J$ is open in $X$, which has relations $P_{15}, av, 2vc$.
Propositions \ref{Y_1 and Y_13 pushforwards} and \ref{Y_2 pushforward} say precisely that
$A_1,A_2,A_3,A_4,A_5,A_6$ are relations on $J$. Therefore each relation on $Q$
is a relation on $J$, and they have the same generators, so there is a surjection
$Q\rightarrow\CH(J)$.

Let $\omega\in Q$ be a homogeneous element of the kernel. By excision,
$\omega$ is pushed forward from $B$. Theorem \ref{degree three} therefore says
that $3\omega$ is a combination of pushforwards from $\mathcal Y_1$, $\mathcal Y_2$, and
$\mathcal Y_{22}$. Their images in $Q$ vanish by Propositions \ref{Y_1 and Y_13 pushforwards},
\ref{Y_2 pushforward}, and \ref{no relations}. Thus $3\omega=0$.

If $\deg\omega\geq2$, then $\omega$ restricts to zero on $J_{\mathrm{nh}}$, so
$\omega\in(v)$ and $112\omega=0$ by Lemma \ref{112 annihilator}.
Hence $\omega$ is killed by both 3 and 112, so $\omega=0$. If $\omega$ is of
degree zero or one, consider
$$
Q^0=\mathbb Z,\qquad Q^1=\mathbb Z\langle a,v\rangle/\langle9a+14v\rangle\cong\mathbb Z.
$$
These are both torsion-free, so $3\omega=0$ again gives $\omega=0$. Therefore the kernel
of $Q\rightarrow\CH(J)$ is trivial and the proposition is proved.
\end{proof}

\begin{proposition}
We have
$$
\CH(H)=\frac{\mathbb Z[b,c,v]}{(2c,14v,2v^2-24b,v^4-35bv^2+24b^2-vc)}.
$$
\end{proposition}

\begin{proof}
All we have to do is excise the non-squarefree locus from
$[B_8/(\PGL_2\times\mathbb G_m)]$, whose Chow ring is
$\mathbb Z[b,c,v]/(2c)$.
By \cite[Corollary 5.3 and the following paragraph]{DL}, with $n=4$,
the image we must excise is generated by pushforwards from $\mathcal Y_1^{\mathrm{sp}}$
and $\mathcal Y_2^{\mathrm{sp}}$ (where Di Lorenzo's $\pi_1$ and $\pi_2$ correspond to our
$f_1$ and $f_2$, respectively). Note that Di Lorenzo deals with the projectivization, and hence
doesn't have to account for the zero section of $B_8$. However, this zero section is in the image of
$\mathcal Y_2^{\mathrm{sp}}$.

We have thus already done all of the relevant work in
Propositions \ref{Y_1 and Y_13 pushforwards} and \ref{Y_2 pushforward}. These give
$$
14v, 2v^2-24b, 120b+60v^2, 50vb-10v^3, v^4-35bv^2+24b^2-vc
$$
as relations, in addition to $2c=0$ coming from $B\PGL_2$. The third and fourth relations aren't
needed, since
$$
\begin{aligned}
120b+60v^2&=5v(14v)-5(2v^2-24b),\\
50vb-10v^3&=-5b(14v)-5v(2v^2-24b).
\end{aligned}
$$
\end{proof}

\begin{proposition}\label{H kernel}
The restriction $\iota^*:\CH(J)\rightarrow\CH(H)$ is surjective, with kernel equal to $(a,2c)$.
\end{proposition}

\begin{proof}
The ring $\CH(H)$ is generated by $b,c,v$, so $\iota^*$ is surjective. In $H$, we have
$a=2c=0$, so the kernel of $\iota^*$ certainly contains $(a,2c)$. We will show that the induced map
$$
\overline{\iota^*}:\CH(J)/(a,2c)\rightarrow\CH(H)
$$
is an isomorphism. Modulo $(a,2c)$, the $A_i$ become
$$
\begin{aligned}
A_1&\mapsto -14v,\qquad &A_2&\mapsto2v^2-24b,\\
A_3&\mapsto 14bv,\qquad &A_4&\mapsto0,\\
A_5&\mapsto 50bv-10v^3,\qquad &A_6&\mapsto v^4-35bv^2+24b^2-vc.
\end{aligned}
$$
Since
$$
14bv=-b(-14v),\qquad
50bv-10v^3=5b(-14v)-5v(2v^2-24b),
$$
the images of the $A_i$ generate precisely the defining ideal of $\CH(H)$ in $\mathbb Z[b,c,v]/(2c)$.
Moreover, $P_{15}$ vanishes modulo $(a,2c)$ by Lemma \ref{ordinary and special Chow rings},
as do $av$ and $2vc$. Thus $\overline{\iota^*}$ is an isomorphism, and so $\ker(\iota^*)=(a,2c)$.
\end{proof}

\section{The Chow rings of $\mathcal A_3^{\mathrm{ind}}$ and $\mg3$}\label{root section}

We now use our computation of $\CH(J)$ to determine the Chow rings of
$\mathcal A_3^{\mathrm{ind}}$ and $\mg3$. We will show that these are a root gerbe over $J$
and a root stack along $H$, respectively. Applying the root gerbe and root stack formulas
then gives their Chow rings, which we express in terms of the Hodge classes.

\subsection{The root stack and the root gerbe}

\begin{lemma}\label{classical blowup}
We have the following identification:
$$
J\cong[(\Bl_ZP)_{\mathrm{good}}/\PGL_3].
$$
\end{lemma}

\begin{proof}
The good locus of $\Bl_ZP$ consists of smooth quartics and squarefree octics over
smooth conics \cite[proof of Proposition 11.1]{CFG}, which is precisely what is left in
our normal cone construction after excision. Locally, the deformation to the normal cone
has equations $f_i=sn_i$ where the $f_i$ generate the ideal of $Z$, and
$t\cdot(s,n_i)=(ts,t^{-1}n_i)$. On the good locus we may normalize by dividing by some
nonzero $n_j$, so that our equations are now $f_i=f_jm_i$, $m_i=n_i/n_j$, which are precisely the
(normalized) equations for the ordinary blowup. Hence
$$
J\cong[(\Bl_ZP)_{\mathrm{good}}/\PGL_3].
$$
\end{proof}

\begin{note}
The preceding lemma raises the question of why we used the deformation to the normal
cone instead of the ordinary blowup. We chose the deformation to the normal cone
purely for organizational reasons, since both arguments require essentially the same excision
calculation. The ordinary and special parts of the deformation to the normal cone have simple
Chow rings, and the top Chern class of the normal bundle of the special part is not a zero-divisor,
so the patching lemma directly computes $\CH(X)$.
This gives a simple ambient presentation and allows us to handle the unwanted loci
uniformly in Sections \ref{envelope section}--\ref{rigidification Chow ring}.
\end{note}

\begin{proposition}\label{A_3 rigidification}
Let $\mathcal A_3^{\mathrm{ind}}\rigidify\mu_2$ be the rigidification of $\mathcal A_3^{\mathrm{ind}}$
along the central inversion. There is an isomorphism
$$
\mathcal A_3^{\mathrm{ind}}\rigidify\mu_2\cong J,
$$
which identifies $\mathcal A_3^{\mathrm{hyp}}\rigidify\mu_2$ with $H$.
\end{proposition}

\begin{proof}
We use the identification from Lemma \ref{classical blowup} to work with the ordinary blowup.
By \cite[proof of Proposition 11.1]{CFG}, the geometric quotient of
$(\Bl_ZP)_{\mathrm{good}}$ by $\PGL_3$ is $M_3$, the coarse moduli space of $\mg3$.
To compare the stack structures we use the characterization of canonical stacks from
\cite[Theorem 1 and Remark 1]{GS}.

The stabilizers of $J$ are $\Aut(C)$ at a non-hyperelliptic curve $C$ and $\Aut(C)/\langle\iota\rangle$
at a hyperelliptic curve with involution $\iota$, by the description of the
quotient of quartics in \cite[\S 7]{DLFV} and of the hyperelliptic quotient in
\cite[Corollary 4.7(b)]{AV}. By \cite[Lemma 3.6]{DLP}, these stabilizers are trivial away
from a closed subset of codimension at least two. The stack $J$ is smooth and of finite type,
since $P$ and $Z$ are smooth. For the linearization of \cite[Proposition 4.1 and its proof]{Artebani},
the good locus is contained in the GIT stable locus: smooth
quartics are stable, and our exceptional divisor consists of squarefree binary octics, which are also
stable \cite[Lemmas 1.4(ii) and 3.4(i)]{Artebani}. Thus the $\PGL_3$-action on
$(\Bl_ZP)_{\mathrm{good}}$ is proper with finite stabilizers, so $J$ is separated and
Deligne--Mumford. Since we work over $\mathbb C$, it is also tame.
The stack $\mathcal A_3^{\mathrm{ind}}\rigidify\mu_2$ is also smooth, separated, and tame with
coarse space $M_3$, and by \cite[\S 4 and Lemma 4.3]{DLP} it has the same stabilizers.
Consequently both $J$ and $\mathcal A_3^{\mathrm{ind}}\rigidify\mu_2$ are the canonical stack over
$M_3$, and \cite[Theorem 1 and Remark 1]{GS} thus identifies them over $M_3$.
This identification carries the reduction of $\mathcal A_3^{\mathrm{hyp}}\rigidify\mu_2$
to the special divisor $H$.
\end{proof}

Let $\rig:\mathcal A_3^{\mathrm{ind}}\rightarrow J$ be the rigidification and let
$\tau:\mg3\rightarrow\mathcal A_3^{\mathrm{ind}}$ be the Torelli morphism.
We will describe $\mathcal A_3^{\mathrm{ind}}$ and $\mg3$ as stacks over $J$, using
root gerbes and root stacks, respectively, which will enable us to compute their Chow rings.

\begin{proposition}\label{A_3^ind root gerbe}
There is an isomorphism
$$
\mathcal A_3^{\mathrm{ind}}\cong J\left(\sqrt{L|_J}\right).
$$
If $M$ denotes the square root of $L|_J$ and $\mathbb E_3$ is the Hodge bundle of
$\mathcal A_3^{\mathrm{ind}}$, then
$$
\mathbb E_3=V^\vee\otimes M^{-1},\qquad
u:=c_1(M),\qquad 2u=v.
$$
\end{proposition}

\begin{proof}
The Torelli morphism $\tau:\mg3\rightarrow\mathcal A_3^{\mathrm{ind}}$ is a finite flat double cover,
simply branched along $\mathcal A_3^{\mathrm{hyp}}$
\cite[\S 4, including the proof of Proposition 4.4]{DLP},
\cite[before Lemma 7.1]{CFG}. Since $\tau$ is a double cover,
\cite[\S 2, Proposition 2.2]{AV} tells us that there is a line bundle $M$ with
$$
\tau_*\mathcal O_{\mg3}=\mathcal O_{\mathcal A_3^{\mathrm{ind}}}\oplus M^{-1},\qquad
M^{\otimes2}=\mathcal O_{\mathcal A_3^{\mathrm{ind}}}(\mathcal A_3^{\mathrm{hyp}})=
\rig^*\mathcal O_J(H)=\rig^*(L|_J).
$$
This defines a morphism
$$
\mathcal A_3^{\mathrm{ind}}\rightarrow J\left(\sqrt{L|_J}\right).
$$
Away from $\mathcal A_3^{\mathrm{hyp}}$, the central involution $[-1]$ exchanges the two sheets of
$\tau$, so it acts as $-1$ on $M^{-1}$ and hence on $M$. Since this character is locally constant,
it must act as $-1$ on $M$ everywhere. Thus the morphism induces an isomorphism on the
relative stabilizers $\mu_2$. Since both stacks are $\mu_2$-gerbes over $J$, it is an
isomorphism, which shows that $\mathcal A_3^{\mathrm{ind}}$ is the root gerbe of $L|_J$.

Now we need to identify the Hodge bundle. By \cite[proof of Lemma 4.5]{DLP}, we know that
$[\mathcal A_3^{\mathrm{hyp}}]=18\lambda_1$. Since $M$ is the square root of
$\mathcal O_{\mathcal A_3^{\mathrm{ind}}}(\mathcal A_3^{\mathrm{hyp}})$,
we have $2u=18\lambda_1$. From our description of $\CH(J)$ in Proposition \ref{excision result},
we know that $\Pic(J)=\langle a,v\rangle/\langle 9a+14v\rangle$. Hence
$\Pic(\mathcal A_3^{\mathrm{ind}})=\langle a,u\rangle/\langle 9a+28u\rangle$. This group is
torsion-free, and so from $2u=18\lambda_1$ we conclude $u=9\lambda_1$ and $a=-28\lambda_1$.
In particular,
\begin{equation}\label{u equality}
M=(\det\mathbb E_3)^9.
\end{equation}

\cite[\S 4, in particular after Corollary 4.1]{vdGK24} constructs a section
$$
\chi_{4,0,8}\in H^0(\mathcal A_3^{\mathrm{ind}},\Sym^4(\mathbb E_3)\otimes(\det\mathbb E_3)^8).
$$
This section is nowhere vanishing by \cite[proof of Proposition 6.1]{CFG}. Recall that
$\mathcal F_4= \Sym^4(V^\vee)\otimes\det V$. Set
$$
W=\mathbb E_3^\vee\otimes(\det\mathbb E_3)^{-9},
$$
so that
$$
\Sym^4(W^\vee)\otimes\det W=\Sym^4(\mathbb E_3)\otimes(\det\mathbb E_3)^8.
$$
The section $\chi_{4,0,8}$ therefore induces a $\GL_3$-equivariant map
$$
\operatorname{Fr}(W)\rightarrow\mathcal F_4\setminus0.
$$
Together with the $\GL_3$-torsor $\operatorname{Fr}(W)\rightarrow\mathcal A_3^{\mathrm{ind}}$,
this defines a morphism
$$
f_\chi:\mathcal A_3^{\mathrm{ind}}\rightarrow[(\mathcal F_4\setminus0)/\GL_3]=[P/\PGL_3].
$$

By \cite[\S 4, in particular after Corollary 4.1]{vdGK24}, $\chi_{4,0,8}$ gives the
equation of the canonical image of the universal curve. Using this description, over the ordinary
locus $f_\chi$ maps to canonical quartics, while over the special locus it maps to doubled
canonical conics (compare with \cite[proof of Proposition 4.2]{Artebani}). That is, $f_\chi$ agrees
with the composition
$$
\mathcal A_3^{\mathrm{ind}}\xrightarrow{\rig}J\rightarrow[P/\PGL_3].
$$
Since we have defined $V$ on $\mathcal A_3^{\mathrm{ind}}$ to mean the pullback of $V$ from
$[P/\PGL_3]$, we see that
$$
V=W=\mathbb E_3^\vee\otimes(\det\mathbb E_3)^{-9}.
$$
Dualizing and rearranging shows that $\mathbb E_3=V^\vee\otimes(\det\mathbb E_3)^{-9}$.
Using \eqref{u equality}, this gives
$$
\mathbb E_3=V^\vee\otimes M^{-1}.
$$
\end{proof}

\begin{proposition}\label{M_3 root stack}
There is an isomorphism over $J$
$$
\mg3\cong J\left(\sqrt H\right)
$$
such that the square root of $\mathcal O_J(H)$ is $\mathcal O_{\mg3}(\mathcal H_3)$.
\end{proposition}

\begin{proof}
Let $f:\mg3\rightarrow J$ be the composite of the Torelli morphism with rigidification.
By Proposition \ref{A_3 rigidification}, the canonical stack of $M_3$ is $J$.

The stack $\mg3$ is smooth, separated and tame, with trivial generic stabilizer.
Its only divisor with nontrivial generic stabilizer is the
hyperelliptic divisor $\mathcal H_3$, whose generic stabilizer is $\mu_2$,
generated by the hyperelliptic involution \cite[Lemma 3.6]{DLP}.
The corresponding divisor on $J$ is $H$. Thus the coarse-space morphism
$\mg3\rightarrow M_3$ has ramification index two along $\mathcal H_3$, and its branch divisor
is the image of $H$ in $M_3$. \cite[Theorem 1]{GS}, applied to $\mg3\rightarrow M_3$,
therefore gives an isomorphism over $J$
$$
\mg3\cong J\left(\sqrt H\right)^{\mathrm{can}},
$$
where $X^{\mathrm{can}}$ is the canonical stack of $X$. Since both $J$ and $H$ are smooth,
the root stack $J\left(\sqrt H\right)$ is smooth and hence is its own canonical stack. We then obtain
$$
\mg3\cong J\left(\sqrt H\right).
$$

Since $f$ has ramification index two along $\mathcal H_3$ over $H$, we have
$$
f^*H=2\mathcal H_3.
$$
Hence
$$
\mathcal O_{\mg3}(\mathcal H_3)^{\otimes 2}\cong f^*\mathcal O_J(H),
$$
and $\mathcal O_{\mg3}(\mathcal H_3)$ is the root line bundle.
\end{proof}

\subsection{Computing the Chow rings}
\begin{proposition}\label{A_ind Chow ring}
We have
$$
\CH(\mathcal A_3^{\mathrm{ind}})=\CH(J)[u]/(2u-v)=
\frac{\mathbb Z[\lambda_1,\lambda_2,\lambda_3]}{I},
$$
where $I$ is the ideal
$$
\left(504\lambda_1^2,24(\lambda_2+10\lambda_1^2),
8\lambda_1(\lambda_2+10\lambda_1^2),28(\lambda_1^3+\lambda_3),
2\lambda_1(2\lambda_1^3+\lambda_1\lambda_2+\lambda_3)\right).
$$
Moreover,
$$
u=9\lambda_1,\qquad a=-28\lambda_1,\qquad b=\lambda_2+261\lambda_1^2,
\qquad c=-\lambda_3-9\lambda_1\lambda_2-810\lambda_1^3.
$$
\end{proposition}

\begin{proof}
The first equality follows from Proposition \ref{excision result},
Proposition \ref{A_3^ind root gerbe} and Lemma \ref{root gerbe Chow}.

By Proposition \ref{A_3^ind root gerbe}, the Hodge bundle is
$$
\mathbb E_3=V^\vee\otimes M^{-1},
$$
where $c_1(M)=u$. Recall that
$$
a=c_1(V),\qquad b=c_2(V),\qquad c=c_3(V).
$$
Let $x_1,x_2,x_3$ be the formal Chern roots of $V$, so that
$$
a=x_1+x_2+x_3,\qquad
b=x_1x_2+x_1x_3+x_2x_3,\qquad
c=x_1x_2x_3.
$$
The Chern roots of $V^\vee$ are $-x_1,-x_2,-x_3$, and tensoring with $M^{-1}$ subtracts
$u$ from each root. Thus the Chern roots of $\mathbb E_3$ are $-x_1-u,-x_2-u,-x_3-u$, and
$$
\begin{aligned}
\lambda_1&=\sum_i(-x_i-u)=-a-3u,\\
\lambda_2&=\sum_{i<j}(-x_i-u)(-x_j-u)=b+2au+3u^2,\\
\lambda_3&=\prod_i(-x_i-u)=-c-ub-u^2a-u^3.
\end{aligned}
$$

After eliminating $v=2u$, the relation $9a+14v=0$ becomes $9a+28u=0$. Using this relation
and the formulas above, \texttt{proposition-7-6.m2} proves the formula for $I$ as well as the
formulas for $a,b,c$.
\end{proof}

\begin{proposition}\label{M_3 Chow ring v1}
Let $h=[\mathcal H_3]\in\CH^1(\mg3)$. Then
$$
\CH(\mg3)=\frac{\CH(J)[h]}{(2h-v, ah, 2ch)}.
$$
\end{proposition}

\begin{proof}
By Proposition \ref{M_3 root stack}, $\mg3\cong J(\sqrt H)$, with root line bundle
$\mathcal O_{\mg3}(\mathcal H_3)$. The restriction $\iota^*:\CH(J)\rightarrow\CH(H)$ is
surjective with kernel $(a,2c)$ by Proposition \ref{H kernel}. Since $[H]=v$ and
$c_1(\mathcal O_{\mg3}(\mathcal H_3))=h$, Lemma \ref{root stack Chow} gives
$$
\CH(\mg3)
=\frac{\CH(J)[h]}{(2h-v,h\ker\iota^*)}
=\frac{\CH(J)[h]}{(2h-v,ah,2ch)}.
$$
\end{proof}

\begin{theorem}\label{M_3 Chow ring}
We have
$$
\CH(\mg3)=\frac{\mathbb Z[\lambda_1,\lambda_2,\lambda_3]}
{(252\lambda_1^2, 24\lambda_2-12\lambda_1^2, 8\lambda_1\lambda_2+80\lambda_1^3,
28\lambda_3+28\lambda_1^3, 2\lambda_1\lambda_3+2\lambda_1^2\lambda_2+4\lambda_1^4)}.
$$
Pullback along the Torelli morphism is surjective and identifies
$$
\CH(\mg3)=\CH(\mathcal A_3^{\mathrm{ind}})/(252\lambda_1^2).
$$
In particular, the order of $\lambda_3$ is 252.
\end{theorem}

\begin{proof}
Plugging the presentation of $\CH(J)$ from Proposition \ref{excision result} into Proposition
\ref{M_3 Chow ring v1} and eliminating $v$ using $v=2h$ gives
$$
\CH(\mg3)=\frac{\mathbb Z[a,b,c,h]}{(P_{15}, ah, 2h c, 9a+28h, 3a^2-24b+8h^2,
ab-28c, 19a^3-88ab+224c)}.
$$
See \texttt{theorem-7-8.m2} for the explicit computation. We now express this in terms of
$\lambda_1,\lambda_2,\lambda_3$.

By Proposition \ref{A_3^ind root gerbe}, $\mathbb E_3=V^\vee\otimes M^{-1}$. Pulling this back
along $\tau$, using $\tau^*M=\mathcal O_{\mg3}(\mathcal H_3)$, gives
$$
\mathbb E_3\cong V^\vee\otimes
\mathcal O_{\mg3}(-\mathcal H_3).
$$
The Chern roots of $\mathbb E_3$ are $-x_1-h,-x_2-h,-x_3-h$, and using precisely the same
computation as in Proposition \ref{A_ind Chow ring} gives
$$
\begin{aligned}
\lambda_1&=-a-3h,\\
\lambda_2&=b+2ah+3h^2,\\
\lambda_3&=-c-hb-h^2a-h^3.
\end{aligned}
$$
The relation $9a+28h=0$ becomes
$$
9(-\lambda_1-3h)+28h=h-9\lambda_1=0.
$$
Eliminating it gives
$$
\begin{aligned}
a&=-28\lambda_1,\qquad h=9\lambda_1,\\
b&=\lambda_2+261\lambda_1^2,\\
c&=-\lambda_3-9\lambda_1\lambda_2-810\lambda_1^3.
\end{aligned}
$$
In particular $ah=-252\lambda_1^2$.

\texttt{theorem-7-8.m2} carries out this change of variables, which proves the claimed presentation,
and computes the order of $\lambda_3$.
Now compare with Proposition \ref{A_ind Chow ring}: the pullback of the universal root $u$ to
$\mg3$ is $h$, and since $2hc=0$ is the pullback of $2uc=0$, which is already a relation on
$\mathcal A_3^{\mathrm{ind}}$,
the only new relation is $ah=-252\lambda_1^2$. This proves the description of the pullback
along the Torelli morphism.
\end{proof}

With rational coefficients the presentation of $\mg3$ reduces to
$$
\CH(\mg3)\otimes\mathbb Q=\mathbb Q[\lambda_1]/(\lambda_1^2),
$$
which agrees with Faber's computation \cite[Theorem 1.9]{Faber}.

Restricting Proposition \ref{M_3 Chow ring v1} to the hyperelliptic locus gives
$$
\CH(\mathcal H_3)=\CH(H)[h]/(2h-v)=\frac{\mathbb Z[b,c,h]}{(2c, 28h, 8h^2-24b)}.
$$
This agrees with the presentation of $\CH(\mathcal H_3)$ for odd $g$ in
\cite[Theorem 6.2]{DL}. The Hodge classes restrict to
$$
\lambda_1=-3h,\qquad \lambda_2=b+3h^2,\qquad
\lambda_3=-c-h b-h^3.
$$

\section{The Chow ring of $A_{1,2}^\circ$}\label{A_12 Chow section}

We now turn to the decomposable locus, beginning with $A_{1,2}^\circ$. We will compute its
Chow ring, normal bundle, and the annihilator of the top Chern class of its normal bundle.
This will show that an element of $\CH(\mathcal A_3\setminus A_{1,1,1})$ is determined
by its restrictions to $\mathcal A_3^{\mathrm{ind}}$ and $A_{1,2}^\circ$, together
with its rationalization.

\begin{definition}\label{A_12 definition}
Let $V_2$ be the standard representation of $\GL_2$ and $\Delta_6$ the locus
of binary sextics with repeated roots. Then
$$
\mg2=\left[(\Sym^6(V_2^\vee)\otimes\det(V_2)^2\setminus\Delta_6)/\GL_2\right],
$$
where $A\in\GL_2$ acts by $(A\cdot f)(X)=\det(A)^2f(A^{-1}X)$, and
$$
\mg{1,1}=[\{(a,b)\in\mathbb A^2:4a^3+27b^2\neq0\}/\mathbb G_m],
$$
where $\mathbb G_m$ acts on $a$ and $b$ with weights $4$ and $6$.

Let $A_{1,2}^\circ=\mg{1,1}\times\mg2$, with closed immersion
$p_{12}:A_{1,2}^\circ\rightarrow\mathcal A_3\setminus A_{1,1,1}$ given by
$(E,C)\mapsto E\times J(C)$.
Let $\mathbb E_1$ and $\mathbb E_2$ be the Hodge bundles of the two factors, and put
$$
t=c_1(\mathbb E_1),\quad\lambda'_1=c_1(\mathbb E_2),\quad \lambda'_2=c_2(\mathbb E_2),\qquad
\alpha_{12}=(p_{12})_*1,\quad \alpha_{12}^{(1)}=(p_{12})_*t.
$$
On $\mathcal A_3\setminus A_{1,1,1}$, $\lambda_1,\lambda_2,\lambda_3$ will still
denote the Chern classes of the Hodge bundle.
\end{definition}

\begin{proposition}\label{A_12 Chow ring}
We have
$$
\CH(A_{1,2}^\circ)=\frac{\mathbb Z[t,\lambda'_1,\lambda'_2]}
{(12t,10\lambda'_1,2{\lambda'_1}^2-24\lambda'_2)},
\qquad N_{A_{1,2}^\circ/(\mathcal A_3\setminus A_{1,1,1})}=\mathbb E_1^\vee\otimes\mathbb E_2^\vee.
$$
If $e$ is the top Chern class of the normal bundle of $A_{1,2}^\circ$ in $\mathcal A_3\setminus A_{1,1,1}$,
then we also have
$$
\begin{aligned}
e&=t^2+t\lambda'_1+\lambda'_2,
& p_{12}^*\lambda_1&=t+\lambda'_1,\\
p_{12}^*\lambda_2&=t\lambda'_1+\lambda'_2,
& p_{12}^*\lambda_3&=t\lambda'_2,\\
p_{12}^*\alpha_{12}&=e,
& p_{12}^*\alpha_{12}^{(1)}&=et.
\end{aligned}
$$
\end{proposition}

\begin{proof}
By \cite{Vistoli}, $\CH(\mg2)=\mathbb Z[\lambda'_1,\lambda'_2]/
(10\lambda'_1,2{\lambda'_1}^2-24\lambda'_2)$.
Then $B\mathbb G_m\times\mg2$ has Chow ring $\CH(\mg2)[t]$, and
$[\mathbb A^2_{4,6}/\mathbb G_m]\times\mg2$ has Chow ring
$\CH(\mg2)[t]$ by vector bundle invariance. We only need to excise the singular locus, which
has weight 12. Therefore
$$
\CH(A_{1,2}^\circ)=\frac{\mathbb Z[t,\lambda'_1,\lambda'_2]}
{(12t, 10\lambda'_1, 2{\lambda'_1}^2-24\lambda'_2)}.
$$
The restrictions of $\alpha_{12}$ and $\alpha_{12}^{(1)}$ follow from the self-intersection formula.

By Lemma \ref{Kodaira-Spencer normal}, the normal bundle of $A_{1,2}^\circ$ is
$\mathbb E_1^\vee\otimes\mathbb E_2^\vee$, from which we get $e=t^2+t\lambda'_1+\lambda'_2$.
The identities for $p_{12}^*\lambda_1$, $p_{12}^*\lambda_2$, $p_{12}^*\lambda_3$
follow from $p_{12}^*\mathbb E_3=\mathbb E_1\oplus\mathbb E_2$.
\end{proof}

\begin{theorem}\label{A_12 annihilator}
We have
$$
\Ann(e)=120\mathbb Z\subseteq\CH^0(A_{1,2}^\circ).
$$
Every element of positive degree of $\CH(A_{1,2}^\circ)$ is killed by 120,
and multiplication by $e$ is injective on $\CH^{>0}(A_{1,2}^\circ)$.
\end{theorem}

\begin{proof}
The statement that all positive degree elements are annihilated by 120 follows from the fact
that $12t=0$, $10\lambda'_1=0$, and $24\lambda'_2=2{\lambda'_1}^2$. The statement that
multiplication by $e$ is injective in positive degree follows from the annihilator calculation,
which is contained in \texttt{theorem-8-3.m2}.
\end{proof}

\begin{lemma}\label{A_12 detection}
A class $\omega\in\CH(\mathcal A_3\setminus A_{1,1,1})$ with
$\omega|_{\mathcal A_3^{\mathrm{ind}}}=0$, $p_{12}^*\omega=0$, and $\omega_{\mathbb Q}=0$ is zero.
\end{lemma}

\begin{proof}
By localization, $\omega=(p_{12})_*z$ for some $z\in\CH(A_{1,2}^\circ)$.
The self-intersection formula gives $ez=p_{12}^*\omega=0$.
By Theorem \ref{A_12 annihilator}, $z=120n$ for some $n\in\mathbb Z$,
so $\omega=120n\alpha_{12}$.
Van der Geer (\cite[Propositions 3.1 and 3.2]{vdG}) gives
$$
\alpha_{12}=\frac{21}{2}\lambda_1^2\neq0
\quad\text{in }\CH^2(\mathcal A_3)_{\mathbb Q}.
$$
Removing the codimension three locus $A_{1,1,1}$ does not affect codimension two,
so this remains nonzero rationally in $\mathcal A_3\setminus A_{1,1,1}$.
Since $\omega_{\mathbb Q}=0$, we have $n=0$ and hence $\omega=0$.
\end{proof}

\section{The Chow rings of $B(\mathbb G_m^3\rtimes S_3)$ and $A_{1,1,1}$}\label{BG_3 section}

We compute $\CH(A_{1,1,1})=\CH([\mg{1,1}^3/S_3])$ using the Weierstrass equation for each of
the three factors of $\mg{1,1}$, with $S_3$ permuting the factors. This leads us to first compute
$\CH(BG_3)$ for $G_3=\mathbb G_m^3\rtimes S_3$. We then excise the locus where at least one
of the three elliptic curves is singular.

\subsection{The Chow ring of $B(\mathbb G_m^3\rtimes S_3)$}
\begin{definition}\label{G_3 definition}
For each of the following groups $G$, denote its standard representation by $\Std_G$.

\begin{enumerate}[label=(\alph*)]
\item Let $T=\mathbb G_m^3$ and $G_3=T\rtimes S_3$. Let the standard representation
be as follows. Let $e_1,e_2,e_3$ be a basis for $V$.
Let $G_3$ act on $V$ as: each copy of $\mathbb G_m$ acts on its respective basis vector,
and $S_3$ permutes. Write
$$
c_j=c_j(\Std_{G_3}).
$$
\item Let
$G_2=\mathbb G_m^2\rtimes S_2=(\{1\}\times\mathbb G_m^2)\rtimes\langle(23)\rangle\leq G_3$.
Define its standard representation similarly to $G_3$. Write
$$
\beta_1=c_1(\Std_{G_2}),\qquad \beta_2=c_2(\Std_{G_2}),\qquad\text{and}\qquad \gamma=c_1(\sgn).
$$
Here $\sgn$ is pulled back from $S_2$. Our $G_2$ and $\Std_{G_2}$ are Larson's
$G$ and $W$, and our $\beta_1,\beta_2,\gamma$ agree with his.
\item Let $H_2=\mathbb G_m^3\rtimes\langle(23)\rangle=\mathbb G_m\times G_2\leq G_3$.
Then $\Std_{H_2}=\Std_{\mathbb G_m}\oplus\Std_{G_2}$. For $H_2$'s first copy of $\mathbb G_m$,
write
$$
w=c_1(\Std_{\mathbb G_m}).
$$
\item For $S_3$, let $\Std_{S_3}$ denote its two-dimensional standard representation,
so that $\Std_{G_3}|_{S_3}=\mathbf 1\oplus\Std_{S_3}$. Write
$$
d_1=c_1(\Std_{S_3})=c_1(\sgn)\qquad\text{and}\qquad d_2=c_2(\Std_{S_3}).
$$
\end{enumerate}
We use the same notation for the pullbacks of these classes along natural maps.
\end{definition}

We will need the Chow ring of $BG_2$, which was computed in \cite{Larson},
and the Chow ring of $BS_3$, which we compute below. For general results on $\CH(BS_3)$,
see \cite[Theorem 11.1]{Totaro}.

\begin{proposition}[{\cite[Theorem 5.2]{Larson}}]\label{G_2 Chow ring}
Using the notation of Definition \ref{G_3 definition}, we have
$$
\CH(BG_2)=\frac{\mathbb Z[\beta_1,\beta_2,\gamma]}{(2\gamma,\gamma(\beta_1+\gamma))}=\frac{\mathbb Z[\bar \beta_1,\beta_2,\gamma]}{(2\gamma,\gamma\bar \beta_1)},
\qquad \bar \beta_1=\beta_1-\gamma.
$$
\end{proposition}

\begin{proposition}
We have
$$
\CH(BS_3)=\mathbb Z[d_1,d_2]/(2d_1,3d_2).
$$
\end{proposition}

\begin{proof}
Let $C_2=\langle(23)\rangle$ and $C_3=\langle(123)\rangle$.
We have
$$
\CH(BC_2)=\mathbb Z[t]/(2t),\qquad
\CH(BC_3)=\mathbb Z[u]/(3u)
$$
where $t=c_1(\sgn)$ and $u=c_1(\chi)$ where $\chi$ maps $(123)$ to a primitive cube root of unity.
The maps $p_2:BC_2\rightarrow BS_3$ and $p_3:BC_3\rightarrow BS_3$ have degrees three
and two, respectively, and hence $(p_2)_*p_2^*x=3x$ and $(p_3)_*p_3^*x=2x$ for all $x$.
Thus any element in $\CH(BS_3)$ which vanishes simultaneously in $\CH(BC_2)$ and
$\CH(BC_3)$ must be annihilated by both 2 and 3, and hence must be trivial.
Thus we have an injection
$$
\rho:\CH(BS_3)\rightarrow\mathbb Z[t]/(2t)\times\mathbb Z[u]/(3u).
$$

The restrictions of $\Std_{S_3}$ are $\mathbf 1\oplus\sgn$ on $C_2$ and $\chi\oplus\chi^{-1}$
on $C_3$. The sign character is trivial on $C_3$, so
$$
\rho(d_1)=(t,0),\qquad \rho(d_2)=(0,-u^2).
$$
Conjugation by a transposition inverts $C_3$, and hence sends $u$ to $-u$.
The restriction of any class from $BS_3$ is invariant under this action,
so its $C_3$ component is a polynomial in $u^2$.
The two restrictions have the same constant term. Therefore every class
has restrictions of the form
$$
\left(n+f(t), n+g(u^2)\right),
$$
where $n\in\mathbb Z$ and $f,g$ are polynomials with no constant term.
This is the restriction of $n+f(d_1)+g(-d_2)$, and so,
by injectivity of $\rho$, the classes $d_1$ and $d_2$ generate
$\CH(BS_3)$.

Finally, the kernel of the composite
$$
\mathbb Z[d_1,d_2]\rightarrow\CH(BS_3)
\xrightarrow{\rho}\mathbb Z[t]/(2t)\times\mathbb Z[u]/(3u)
$$
is
$$
(d_2,2d_1)\cap(d_1,3d_2)=(2d_1,3d_2).
$$
Since $\rho$ is injective, these are exactly the relations among
$d_1$ and $d_2$.
\end{proof}

Consider the stratification of $\Std_{G_3}$ by the number of nonzero coordinates.
Let $V_{\geq k}\subseteq \Std_{G_3}$ be the open subset with at least $k$ nonzero coordinates,
and put $\mathcal V_k=[V_{\geq k}/G_3]$. Then
$$
\mathcal V_1=\mathcal V_3\sqcup\left(\mathcal V_2\setminus\mathcal V_3\right)
\sqcup\left(\mathcal V_1\setminus\mathcal V_2\right).
$$
The group $G_3$ acts transitively on each stratum, so it suffices
to check the stabilizers at specific points.
The stabilizers of $(1,1,1)$, $(0,1,1)$, $(1,0,0)$, and $(0,0,0)$ are, respectively, $S_3$,
$\mathbb G_m\times S_2$, $G_2$, and $G_3$. Hence
$$
\mathcal V_3=BS_3,\qquad
\mathcal V_2\setminus\mathcal V_3=B(\mathbb G_m\times S_2),\qquad
\mathcal V_1\setminus\mathcal V_2=BG_2.
$$
The top Chern classes of the normal bundles are as follows:
\begin{itemize}
\item $\mathcal V_2\setminus\mathcal V_3$: the normal bundle at $(0,1,1)$ is spanned by
the first coordinate. The stabilizer $\mathbb G_m\times S_2$ of $(0,1,1)$ acts on the normal bundle
with weight one from $\mathbb G_m$ and trivially from $S_2$. That is, the normal bundle
is the standard representation of $\mathbb G_m$, whose top Chern class is
$w\in\CH(B(\mathbb G_m\times S_2))=\mathbb Z[\gamma,w]/(2\gamma)$, which is not a zero-divisor.
\item $\mathcal V_1\setminus\mathcal V_2$: the normal bundle at $(1,0,0)$ is spanned by
the second and third coordinates. The stabilizer $G_2$ of $(1,0,0)$ acts on these as its standard
representation from \cite[Definition 2.3]{Larson}, whose top Chern class is
$\beta_2\in\CH(BG_2)=\mathbb Z[\bar \beta_1,\beta_2,\gamma]/(2\gamma,\gamma\bar \beta_1)$,
which is not a zero-divisor.
\end{itemize}

We will consider
$$
[\mathbb P(\Std_{G_3})/G_3]=[V_{\geq1}/(G_3\times\mathbb G_m)]
$$
where the action is $(g,t)\cdot v=tg(v)$. Under the change of coordinates $(g,t)\mapsto(t^{-1}g,t)$
we see that $\mathbb G_m$ acts trivially. Hence
\begin{equation}\label{P direct product}
[\mathbb P(\Std_{G_3})/G_3]=\mathcal V_1\times B\mathbb G_m.
\end{equation}
We will compute the Chow ring of $\mathcal V_1\times B\mathbb G_m$ then
unprojectivize to find the Chow ring of $BG_3$.

\begin{proposition}\label{V_1 Chow ring}
The Chow ring of $\mathcal V_1$ is
$$
\CH(\mathcal V_1)=\frac{\mathbb Z[u_1,u_2,d_1,d_2]}
{(2d_1,3d_2,d_2 u_1,d_2 u_2,d_1 u_1u_2)}.
$$
Here $u_1$ is the class of the closure of $\mathcal V_2\setminus\mathcal V_3$ in $\mathcal V_1$
and $u_2$ is the class of $\mathcal V_1\setminus\mathcal V_2$. That is, they are the classes of the
union of the coordinate planes and the union of the coordinate axes, respectively.
\end{proposition}

\begin{proof}
First we patch $\mathcal V_2\setminus\mathcal V_3$ into $\mathcal V_3$. By the
discussion before this proposition, the normal bundle of $\mathcal V_2\setminus\mathcal V_3$
has top Chern class $w$, which is not a zero-divisor. Hence by Lemma \ref{patching lemma},
we have that $\CH(\mathcal V_2)$ is the following fiber product:
$$
\begin{tikzcd}[column sep=large, row sep=large]
\displaystyle\CH(\mathcal V_2)\arrow[r]\arrow[d]&
\displaystyle\frac{\mathbb Z[\gamma,w]}{(2\gamma)}\arrow[d,"w\mapsto 0"]\\
\displaystyle\frac{\mathbb Z[d_1,d_2]}{(2d_1,3d_2)}
\arrow[r,"\substack{d_1\mapsto \gamma\\ d_2\mapsto 0}"']&
\displaystyle\frac{\mathbb Z[\gamma]}{(2\gamma)}.
\end{tikzcd}
$$
Put
$$
R=\mathbb Z[u_1,d_1,d_2],\qquad
I=(u_1,2d_1,3d_2),\qquad J=(d_2,2d_1).
$$
Since $(2d_1,3d_2)\subseteq J$ and $u_1$ isn't a zero-divisor in $R/J$, Lemma
\ref{fiber product lemma} gives
$$
I\cap J=(2d_1,3d_2)+u_1J=(2d_1,3d_2,d_2u_1).
$$
Also $I+J=(u_1,d_2,2d_1)$, so the same lemma gives
$$
\CH(\mathcal V_2)=\frac{\mathbb Z[u_1,d_1,d_2]}{(2d_1,3d_2,d_2u_1)}.
$$
Here
$$
d_1=(d_1,\gamma)\qquad\text{and}\qquad d_2=(d_2,0)
$$
are the pullbacks of $d_1$ and $d_2$ from $BS_3$ to $\mathcal V_2$ along the map
$$
\mathcal V_2\rightarrow[\Std_{G_3}/G_3]\rightarrow BG_3\rightarrow BS_3,
$$
and
$$
u_1=(0,w)
$$
is the fundamental class of the divisor $\mathcal V_2\setminus\mathcal V_3$.

We can find expressions for the restrictions of $c_1$ and $c_2$ to $\mathcal V_2$, as follows.
If $x_1,x_2,x_3$ are the coordinate functions, the divisor $\mathcal V_2\setminus\mathcal V_3$
is precisely the vanishing of $x_1x_2x_3$. Hence its fundamental class, $u_1$, is equal to
$c_1(\chi)$, where $\chi$ is the character $\chi(t_1,t_2,t_3,\sigma)=t_1t_2t_3$. Since
$$
\det\Std_{G_3}=\chi\otimes\sgn,
$$
we see that $c_1=u_1+d_1$ on $\mathcal V_2$. For $c_2$, we compute the
restriction of $\Std_{G_3}$ to $\mathcal V_3=BS_3$ and
$\mathcal V_2\setminus\mathcal V_3=B\mathbb G_m\times BS_2$. Its restriction to
$BS_3$ is just the permutation representation $P_3$, and its restriction to
$B\mathbb G_m\times BS_2$ is $\Std_{\mathbb G_m}\oplus\mathbf 1\oplus\sgn$. Thus, in the
fiber product, the second Chern class is $(d_2,\gamma w)$. Expressed with our chosen generators,
this says
$$
c_2=d_1 u_1+d_2\quad\text{in }\CH(\mathcal V_2).
$$

We now patch in $\mathcal V_1\setminus\mathcal V_2$. By the discussion before this
proposition, the normal bundle has top Chern class $\beta_2$, which is not a zero-divisor.
On $\mathcal V_1\setminus\mathcal V_2$, we have $\Std_{G_3}=\mathbf 1\oplus\Std_{G_2}$.
Therefore $u_1$ restricts to $\bar \beta_1$, $d_1$ to $\gamma$, and $d_2$ to zero.
Hence by the patching lemma we have that $\CH(\mathcal V_1)$ is the following fiber product:
$$
\begin{tikzcd}[column sep=large, row sep=large]
\displaystyle
\CH(\mathcal V_1)\arrow[r]\arrow[d]
&
\displaystyle\frac{\mathbb Z[\bar \beta_1,\beta_2,\gamma]}{(2\gamma, \gamma\bar \beta_1)}
\arrow[d,"\beta_2\mapsto 0"]\\
\displaystyle\frac{\mathbb Z[u_1,d_1,d_2]}{(2d_1,3d_2,d_2 u_1)}
\arrow[r,"\substack{d_1\mapsto \gamma\\ d_2\mapsto 0\\ u_1\mapsto\bar \beta_1}"']
&
\displaystyle\frac{\mathbb Z[\bar \beta_1, \gamma]}{(2\gamma, \gamma\bar \beta_1)}.
\end{tikzcd}
$$
Put
$$
R=\mathbb Z[u_1,u_2,d_1,d_2],\qquad
I=(u_2,2d_1,3d_2,d_2u_1),\qquad J=(d_2,2d_1,d_1u_1).
$$
Since $(2d_1,3d_2,d_2u_1)\subseteq J$ and $u_2$ isn't a zero-divisor in $R/J$,
Lemma \ref{fiber product lemma} gives
$$
\begin{aligned}
I\cap J&=(2d_1,3d_2,d_2u_1)+u_2J\\
&=(2d_1,3d_2,d_2u_1,d_2u_2,d_1u_1u_2).
\end{aligned}
$$
Also $I+J=(u_2,d_2,2d_1,d_1u_1)$, so the same lemma gives
$$
\CH(\mathcal V_1)=\frac{\mathbb Z[u_1,u_2,d_1,d_2]}
{(2d_1,3d_2,d_2u_1,d_2u_2,d_1u_1u_2)}.
$$
Here $u_2=(0,\beta_2)$ is the fundamental class of $\mathcal V_1\setminus\mathcal V_2$.
Comparing both restrictions gives
$$
c_2=u_2+d_1 u_1+d_2\quad\text{in }\CH(\mathcal V_1).
$$
\end{proof}

\begin{proposition}\label{BG Chow ring}
Using the notation of Definition \ref{G_3 definition}, there is an isomorphism
$$
\CH(BG_3)=
\frac{\mathbb Z[c_1,c_2,c_3,d_1,d_2]}{
(2d_1, 3d_2, d_1d_2, d_2 c_1, d_2(c_2-d_2),
d_1[c_3+c_1c_2+d_1(c_1^2+c_2)+d_1^3])}.
$$
Under $\CH([V/G_3])=\CH(BG_3)$, the classes
$$
u_1=c_1-d_1,\qquad u_2=c_2-d_1 u_1-d_2
$$
are, respectively, the classes of the union of the coordinate planes and the union of the coordinate
axes.
\end{proposition}

\begin{proof}
Write
$$
\mathbf P=[\mathbb P(\Std_{G_3})/G_3].
$$
Due to our earlier coordinate change,
$(g,t)\mapsto(t^{-1}g,t)$, the map to $BG_3$ is induced by $(g,t)\mapsto t^{-1}g$.
Let $\psi:\mathbf P\rightarrow BG_3$ be this map.

Using \eqref{P direct product},
$$
\CH(\mathbf P)=\CH(\mathcal V_1\times B\mathbb G_m)=
\CH(\mathcal V_1)[\xi].
$$
Under this equivalence, $\xi=c_1(\mathcal O_{\mathbf P}(1))$ and $\psi^*\Std_{G_3}$ is
the standard bundle on $\mathcal V_1$ tensored with $\mathcal O_{\mathbf P}(-1)$.
Hence this ring is generated by $d_1,d_2,c_1,c_2,\xi$.
Using our computations of $c_1$ and $c_2$ on $\mathcal V_1$ and applying
a twist by $\mathcal O_{\mathbf P}(-1)$ gives, on $\mathbf P$,
$$
\begin{aligned}
c_1&=u_1+d_1-3\xi,\\
c_2&=u_2+d_1 u_1+d_2-2\xi(u_1+d_1)+3\xi^2.
\end{aligned}
$$

The projective bundle formula gives us a second way to compute $\CH(\mathbf P)$:
$$
\CH(\mathbf P)=\CH(BG_3)[\xi]/(\xi^3+c_1\xi^2+c_2\xi+c_3).
$$
Now let $R$ be the subring of $\CH(BG_3)$ generated by $d_1,d_2,c_1,c_2,c_3$.
Then $\CH(\mathbf P)$ is generated by $R$ and $\xi$, and so $\CH(BG_3)$, being
the part of $\CH(\mathbf P)$ which is degree zero in the variable $\xi$ (after expressing
all elements in the form $a_0+a_1\xi+a_2\xi^2$), is generated by $d_1,d_2,c_1,c_2,c_3$.

We now only need to determine the relations among these generators. That is,
we only need to determine the kernel of the map
$$
\mathbb Z[d_1,d_2,c_1,c_2,c_3]\rightarrow\CH(\mathbf P).
$$
This is done in \texttt{proposition-9-5.m2}.
\end{proof}

\subsection{The Chow ring of $A_{1,1,1}$}

\begin{definition}\label{A_1,1,1 quotient}
Let $V_m$ be the representation $\Std_{G_3}$, but with each $\mathbb G_m$ acting with weight $m$.
Then, giving $V_4\oplus V_6$ the coordinates $(a_i,b_i)$, we have
$$
A_{1,1,1}=\left[\frac{\{(a_i,b_i):4a_i^3+27b_i^2\neq0\}}{G_3}\right].
$$
For a subset $I\subseteq\{1,2,3\}$ of cardinality $r$, let $G_I$ be its setwise stabilizer in $G_3$
and put
$$
\mathcal N_I=[(\mathbb A^1_2)^I\times(\mathbb A^2_{4,6})^{I^c}/G_I].
$$
Note that $G_{\{1\}}=\mathbb G_m\times(\mathbb G_m^2\rtimes S_2)=H_2$.
There is a proper map
$$
\mathcal N_I\rightarrow[(V_4\oplus V_6)/G_3]
$$
given by $v_i\mapsto(-3v_i^2,2v_i^3)$ for $i\in I$, and by the identity on the remaining factors.
Up to permutation, there is one such $\mathcal N_r$ for each $r=1,2,3$. We take
$$
\mathcal N_1=\mathcal N_{\{1\}},\qquad\mathcal N_2=\mathcal N_{\{2,3\}},\qquad
\mathcal N_3=\mathcal N_{\{1,2,3\}}.
$$
\end{definition}

\begin{proposition}\label{A_1,1,1 Chow ring}
We have
$$
\CH(A_{1,1,1})=\CH(BG_3)/(12c_1,24c_2,36c_3),
$$
and additively
$$
\CH(A_{1,1,1})=\frac{\mathbb Z[c_1,c_2,c_3]}{(12c_1,24c_2,36c_3)}
\oplus d_1\mathbb F_2[d_1,c_1,c_2]\oplus d_2\mathbb F_3[d_2,c_3].
$$
\end{proposition}

\begin{proof}
We need to compute the image of the boundary in $\CH([(V_4\oplus V_6)/G_3])=\CH(BG_3)$.
The $\mathcal N_r$ of Definition \ref{A_1,1,1 quotient} form a Chow envelope for the boundary
since the normalization map for the cuspidal cubic $4a_i^3+27b_i^2=0$ does.

Consider the map
$$
\pi: BH_2\rightarrow BG_3.
$$
There's an induced map on vector bundles
$$
\mathcal N_1\rightarrow[(V_4\oplus V_6)/H_2]\rightarrow[(V_4\oplus V_6)/G_3].
$$
The discriminant of the first factor has weight 12, so its divisor has class $12w$
(for the definition of $w$, see Definition \ref{G_3 definition}).
Since $\mathcal N_1$ is a vector bundle over $BH_2$, the projection formula says the image
of $\CH(\mathcal N_1)$ in $\CH([(V_4\oplus V_6)/H_2])$ is $(12w)$.
Therefore the image in $\CH([(V_4\oplus V_6)/G_3])$ is
$$
\mathfrak J_1=\pi_*\left((12w)\right).
$$
On $H_2$ we have $\pi^*\Std_{G_3}=\Std_{\mathbb G_m}\oplus\Std_{G_2}$,
and hence
$$
\pi^*c_1=\beta_1+w,\qquad \pi^*c_2=\beta_2+\beta_1w,\qquad \pi^*c_3=\beta_2w,\qquad
\gamma=\pi^*d_1.
$$
This allows us to express $\beta_1$ and $\beta_2$ in terms of $w$:
\begin{equation}\label{helpful formula}
\beta_1=\pi^*c_1-w,\qquad \beta_2=\pi^*c_2-\pi^*c_1w+w^2.
\end{equation}
Substituting this expression for $\beta_2$ into $\pi^*c_3=\beta_2w$ gives
$$
w^3-\pi^*c_1w^2+\pi^*c_2w-\pi^*c_3=0,
$$
and so $\CH(BH_2)$ is generated by $1,w,w^2$ as a $\CH(BG_3)$-module.
So to compute the pushforward of $12w\cdot a\in(12w)$, it suffices to compute for
$a=1,w,w^2$. Hence $\mathfrak J_1$ is generated by $\pi_*(12w)$, $\pi_*(12w^2)$ and
$\pi_*(12w^3)$.

Now consider the cover
$$
\pi':B\mathbb G_m^3\rightarrow BG_3.
$$
On $B\mathbb G_m^3$, write $x_1,x_2,x_3$ for the generators of its Chow ring.
Then
$$
\pi'^*c_1=x_1+x_2+x_3,\qquad \pi'^*c_2=x_1x_2+x_1x_3+x_2x_3,\qquad \pi'^*c_3=x_1x_2x_3.
$$
By the description of $\CH(BG_3)$ from Proposition \ref{BG Chow ring},
the kernel of restriction from $\CH(BG_3)$ to $\CH(B\mathbb G_m^3)$ is killed by 6.
Since we need only compute $\pi_*(12w),\pi_*(12w^2),\pi_*(12w^3)$, each of which
is a multiple of 12, we see that it suffices to compute the restriction of these classes to
$B\mathbb G_m^3$.

Pushing from $BH_2$ to $BG_3$ and then pulling to $B\mathbb G_m^3$ sums over the three
possible preimages, and so $\pi'^*\pi_*(w^k)=x_1^k+x_2^k+x_3^k$.
Thus
$$
\begin{aligned}
\pi_*(12w)&=12(x_1+x_2+x_3)=12c_1,\\
\pi_*(12w^2)&=12(x_1^2+x_2^2+x_3^2)=12c_1^2-24c_2,\\
\pi_*(12w^3)&=12(x_1^3+x_2^3+x_3^3)=12c_1^3-36c_1c_2+36c_3.
\end{aligned}
$$
Consequently $\mathfrak J_1=(12c_1,24c_2,36c_3)$.

This is just the image of $\CH(\mathcal N_1)$, but we claim that the remaining $\mathcal N_I$
add no relations. Consider $r=2$. On $[(V_4\oplus V_6)/H_2]$, the tuple
$(4a_2^3+27b_2^2,4a_3^3+27b_3^2)$ is acted on by $\mathbb G_m^2$ with weight 12
and by $(23)\in H_2$ by swapping. This is a section of the standard representation of $G_2$
with weight 12, pulled back to $BH_2$:
$$
(t_1,t_2,t_3,\sigma)\cdot(z_2,z_3)=(t_2^{12}\sigma z_2,t_3^{12}\sigma z_3).
$$
The image of $\mathcal N_2$ is the locus where both entries of this tuple
vanish, so to show that $\mathcal N_2$ contributes no further relations, we need to
understand the class of the zero section of this weight-12 representation.

The map from $\Std_{G_2}$ to this representation given by
$$
f:(z_2,z_3)\longmapsto(z_2^{12},z_3^{12})
$$
is equivariant and finite flat. Its inverse image of the zero section is supported on the zero section
with multiplicity $12^2$. Using vector bundle invariance, $f^*$ is the identity, and so we see that the
zero section of the weight-12 representation has class divisible by $12^2$. The projection formula
therefore shows that every pushforward from $\mathcal N_2$ is divisible by $12^2$ in $\CH(BG_3)$.
An analogous argument shows that every class pushed from $\mathcal N_3$ is divisible by $12^3$.

Since $d_1$ and $d_2$ are annihilated by 2 and 3, respectively, any pushforward from
$\mathcal N_2$ or $\mathcal N_3$ must be a polynomial in $c_1,c_2,c_3$. Since $12^r$, $r=2,3$,
is divisible by $12,24,36$, we conclude every pushforward lies in
$\mathfrak J_1=(12c_1,24c_2,36c_3)$. This proves the presentation of $\CH(A_{1,1,1})$.

With respect to the lexicographic order $c_3>c_2>c_1>d_1>d_2$, the computation in
\texttt{proposition-9-7.m2} gives the leading-term ideal
$$
(12c_1,24c_2,36c_3,2d_1,3d_2,d_1d_2,d_1c_3,d_2c_1,d_2c_2).
$$
Since the first five generators are themselves relations, reduction by this
Gr\"obner basis gives the displayed additive decomposition.
\end{proof}

\section{The normal bundle of $A_{1,1,1}$}\label{A_1,1,1 normal bundle section}

We next compute the normal bundle of $A_{1,1,1}$, the annihilator of its top Chern class,
and the restrictions of $\alpha_{12}$ and $\alpha_{12}^{(1)}$. These computations will be used to
find generators for $\CH(\mathcal A_3)$ and to prove the injectivity statement in Section
\ref{final ring section}.

\begin{definition}\label{pi definition}
For $p_{111}:A_{1,1,1}\rightarrow\mathcal A_3$ the inclusion, let $\alpha_{111}=(p_{111})_*1$ and $\alpha_{111}^{(1)}=(p_{111})_*d_1$.
Let $A'=[\mathcal A_1^3/S_2]=\mathcal A_1\times[\mathcal A_1^2/S_2]$, where the first
curve is distinguished (i.e. $S_2$ acts by swapping the second and third curves).
Consider the map $A'\rightarrow[\mathcal A_1^3/S_3]$ which forgets which factor
is distinguished. Using the notation of the previous section, this is the same as the map
$$
\left[\frac{\{(a_i,b_i):4a_i^3+27b_i^2\neq 0\}}{H_2}\right]\rightarrow
\left[\frac{\{(a_i,b_i):4a_i^3+27b_i^2\neq 0\}}{G_3}\right]
$$
induced by $\pi:BH_2\rightarrow BG_3$, and so we also refer to this map as $\pi$.
Let $L_i$ be the Hodge bundle of the $i^{\text{th}}$ component of $\mathcal A_1^3$.
Note that $L_1\otimes L_2\otimes L_3$ descends to $A_{1,1,1}$
since it is $S_3$-equivariant. We call its descent $L_{123}$.
The line bundle $L_1$ descends to $A'$, and we call its descent $L_1$ and
its first Chern class $w$.
\end{definition}

\begin{proposition}\label{A_1,1,1 restrictions}
The normal bundle of $A_{1,1,1}$ is
$N_{A_{1,1,1}/\mathcal A_3}=p_{111}^*\mathbb E_3\otimes L_{123}^{-1}$,
and its top Chern class is
$$
e=c_3-c_1c_2+d_1(c_1^2+c_2)+d_1^3.
$$
The pushforwards along $\pi$ of $1,w,w^2$ are
\begin{equation}\label{w powers pushforward}
\pi_*1=3,\quad \pi_*w=c_1+d_1,\quad
\pi_*w^2=c_1^2-2c_2+d_1^2+2d_2.
\end{equation}
The pullbacks of $\alpha_{12}$ and $\alpha_{12}^{(1)}$ to $A_{1,1,1}$ are
\begin{equation}\label{i^* equalities}
\begin{aligned}
p_{111}^*\alpha_{12}&=c_1^2+c_2+d_1^2+2d_2,\\
p_{111}^*\alpha_{12}^{(1)}&=c_1^3-2c_1c_2+3c_3+d_1^2c_1.
\end{aligned}
\end{equation}
\end{proposition}

\begin{proof}
Note that the lift of $p_{111}^*\mathbb E_3$ to $\mathcal A_1^3$ is $L_1\oplus L_2\oplus L_3$,
which descends to $A_{1,1,1}$. Abusing notation, we will call its descent
$L_1\oplus L_2\oplus L_3$. Hence
$$
p_{111}^*\mathbb E_3=L_1\oplus L_2\oplus L_3,\qquad\det p_{111}^*\mathbb E_3=L_{123}\otimes\sgn.
$$
Then
$$
c_1=c_1(\det p_{111}^*\mathbb E_3)=c_1(L_{123})+c_1(\sgn)=c_1(L_{123})+d_1,
$$
and so, using $2d_1=0$ in $\CH(A_{1,1,1})$,
$$
c_1(L_{123})=c_1+d_1.
$$
By Lemma \ref{Kodaira-Spencer normal}, the normal bundle of $A_{1,1,1}$ lifted to
$\mathcal A_1^3$ is
\begin{equation}\label{A_1,1,1 normal lift}
\bigoplus_{i<j}L_i^{-1}\otimes L_j^{-1}=\bigoplus_i L_i\otimes(L_1^{-1}\otimes L_2^{-1}
\otimes L_3^{-1}),
\end{equation}
which descends to $A_{1,1,1}$ as $p_{111}^*\mathbb E_3\otimes L_{123}^{-1}$.
Hence
$$
N_{A_{1,1,1}/\mathcal A_3}=p_{111}^*\mathbb E_3\otimes L_{123}^{-1}.
$$
The top Chern class is
$$
e=c_3-c_1c_2+d_1(c_1^2+c_2)+d_1^3,
$$
which follows by using $2d_1=0$ in
$$
c_3(p_{111}^*\mathbb E_3\otimes L_{123}^{-1})=c_3-(c_1+d_1) c_2+(c_1+d_1)^2c_1-(c_1+d_1)^3.
$$
In our new context, the last relation in Proposition \ref{BG Chow ring} says
that $d_1 e=0$. Thus $e$, unfortunately, is a zero-divisor.

Since $\pi$ is degree three, we have
$$
\pi_*1=3.
$$
By the norm construction and its compatibility with pushforward
\cite[Tags \href{https://stacks.math.columbia.edu/tag/0BCY}{0BCY},
\href{https://stacks.math.columbia.edu/tag/0BD2}{0BD2}, and
\href{https://stacks.math.columbia.edu/tag/02RT}{02RT}]{Stacks}, we have
$$
\begin{aligned}
\pi_*w&=\pi_*c_1(L_1)\\
&=c_1\left(\Nm_{\pi}(L_1)\right)\\
&=c_1\left(\det(\pi_*L_1)\otimes\det(\pi_*\mathcal O_{A'})^{-1}\right)\\
&=c_1(\pi_*L_1)-c_1(\pi_*\mathcal O_{A'})\\
&=c_1(p_{111}^*\mathbb E_3)-c_1(\sgn)\\
&=c_1-d_1\\
&=c_1+d_1.
\end{aligned}
$$
We claim that $\pi_*\beta_2=u_2$. To see this, we will compute on $[V/H_2]\rightarrow [V/G_3]$.
Note that $\beta_2$ is represented by the vanishing of the
coordinates $x_2,x_3$ in $V$. Under $[V/H_2]\rightarrow[V/G_3]$, this maps isomorphically
(away from the origin) onto the union of the coordinate axes, which is the class that represents
$u_2$ by Proposition \ref{BG Chow ring}.

We can now use $\beta_2=\pi^*c_2-\pi^*c_1w+w^2$ from \eqref{helpful formula} to obtain
$w^2=\beta_2-\pi^*c_2+\pi^*c_1w$, and hence
$$
\pi_*w^2=\pi_*\beta_2-\pi_*\pi^*c_2+\pi_*(w\pi^*c_1)=u_2-3c_2+c_1\pi_*w.
$$
Plugging in $u_1=c_1-d_1$, $u_2=c_2-d_1 u_1-d_2$, and
$\pi_*w=c_1+d_1$ proves the claimed formula.

Consider the Cartesian square
$$
\begin{tikzcd}[column sep=large, row sep=large]
A'\arrow[r]\arrow[d,"\pi"'] & \mathcal A_1\times\mathcal A_2\arrow[d,"p_{12}"]\\
A_{1,1,1}\arrow[r,"p_{111}"'] & \mathcal A_3
\end{tikzcd}
$$
where the top horizontal map is the identity on the first factor and the product morphism on the
other two factors. We compute the excess bundle for this square.
By Lemma \ref{Kodaira-Spencer normal}, the lift of $N_{A'/\mathcal A_1\times\mathcal A_2}$
to $\mathcal A_1^3$ is
$$
N_{A'/(\mathcal A_1\times\mathcal A_2)}=L_2^{-1}\otimes L_3^{-1}.
$$
This is $S_2$-equivariant and hence descends to $A'$. Using \eqref{A_1,1,1 normal lift},
the excess bundle $\mathcal E$ therefore fits into the exact sequence
$$
0\rightarrow L_2^{-1}\otimes L_3^{-1}\rightarrow\bigoplus_{1\leq i<j\leq3}L_i^{-1}\otimes L_j^{-1}
\rightarrow\mathcal E\rightarrow0.
$$
Writing $\mathcal H$ for the bundle which descends from $L_2\oplus L_3$, this gives
$$
\mathcal E= L_1^{-1}\otimes\mathcal H^\vee,\qquad
c_2(\mathcal E)=w^2+w c_1(\mathcal H)+c_2(\mathcal H).
$$
Since
$$
\pi^*p_{111}^*\mathbb E_3=L_1\oplus L_2\oplus L_3=L_1\oplus\mathcal H,
$$
we have
$$
\pi^*c_2=c_2(L_1\oplus\mathcal H)
=w c_1(\mathcal H)+c_2(\mathcal H).
$$
Hence
$$
c_2(\mathcal E)=w^2+\pi^*c_2.
$$

Therefore the excess intersection formula gives
$$
\begin{aligned}
p_{111}^*\alpha_{12}&=\pi_*c_2(\mathcal E)=\pi_*(w^2+\pi^*c_2)=\pi_*w^2+3c_2,\\
p_{111}^*\alpha_{12}^{(1)}&=\pi_*(w c_2(\mathcal E))=\pi_*w^3+c_2\pi_*w.
\end{aligned}
$$
Since $\mathcal H$ has rank two,
$$
0=c_3(\mathcal H)
=\pi^*c_3-w\pi^*c_2+w^2\pi^*c_1-w^3.
$$
Solving this for $w^3$ and pushing forward gives
$$
\pi_*w^3=3c_3-c_2\pi_*w+c_1\pi_*w^2.
$$
We can now substitute our values for $\pi_*w$ and $\pi_*w^2$ from \eqref{w powers pushforward}
and use $d_2 c_1=0$ from Proposition \ref{BG Chow ring}, to obtain
$$
\begin{aligned}
p_{111}^*\alpha_{12}&=c_1^2+c_2+d_1^2+2d_2,\\
p_{111}^*\alpha_{12}^{(1)}&=c_1^3-2c_1c_2+3c_3+d_1^2c_1.
\end{aligned}
$$
\end{proof}

\begin{lemma}\label{A_1,1,1 annihilator}
The annihilator of the top Chern class of the normal bundle of $A_{1,1,1}$ in $\CH(A_{1,1,1})$ is the
ideal $\Ann(e)=(36,d_1)$. Its positive-degree part is killed by $2$.
\end{lemma}

\begin{proof}
See \texttt{lemma-10-3.m2}.
\end{proof}

\section{$\mg2^{\mathrm{ct}}$ as a test family}\label{M_2^ct test family section}

Unfortunately, the annihilator of the top Chern class of the normal bundle of $A_{1,1,1}$ contains
torsion, so the argument of Lemma \ref{A_12 detection} does not give injectivity.
To close this gap, we use $\mg2^{\mathrm{ct}}$ as a test family, mapping it to $\mathcal A_3$ by
sending the curve $C$ to the product $E_0\times J(C)$ for a fixed smooth elliptic curve $E_0$.

\begin{definition}
Let $\mg2^{\mathrm{ct}}$ be the stack of stable curves of genus two of compact type,
and let $\Delta_1=[\mg{1,1}^2/S_2]\subseteq\mg2^{\mathrm{ct}}$ be the locus of curves with a
separating node, with inclusion $j:\Delta_1\rightarrow \mg2^{\mathrm{ct}}$.
As above, we write $\lambda'_1,\lambda'_2$ for the first and second Hodge classes of
$\mg2^{\mathrm{ct}}$. These are the $\lambda_i$'s of \cite{HJL}, and we use $'$ to distinguish them
from the $\lambda_i$'s on $\mathcal A_3$. We write $\delta_1=[\Delta_1]$.

Fix a smooth elliptic curve $E_0$ and let $f:\mg2^{\mathrm{ct}}\rightarrow\mathcal A_3$
be the morphism $C\mapsto E_0\times J(C)$.
\end{definition}

\begin{proposition}\label{M_2^ct Chow ring prop}
The Chow ring of $\mg2^{\mathrm{ct}}$ is
$$
\CH(\mg2^{\mathrm{ct}})=\frac{\mathbb Z[\lambda'_1,\lambda'_2,\delta_1]}
{(24{\lambda'_1}^2-48\lambda'_2, 24\delta_1\lambda'_2, 10\lambda'_1-2\delta_1,
2{\lambda'_1}^2-24\lambda'_2+\delta_1\lambda'_1-\delta_1^2)},
$$
and the sign class $\gamma$ on $\Delta_1=[\mg{1,1}^2/S_2]$ satisfies
$$
\gamma=j^*(\delta_1+\lambda'_1),\qquad
j_*\gamma=\delta_1(\delta_1+\lambda'_1)=12({\lambda'_1}^2-2\lambda'_2).
$$
\end{proposition}

\begin{proof}
The Chow ring follows from \cite[Theorem 5.6]{HJL} after changing the first generator using
$$
24{\lambda'_1}^2-48\lambda'_2=12\delta_1\lambda'_1+(2\lambda'_1-\delta_1)
(10\lambda'_1-2\delta_1)+2(2{\lambda'_1}^2-24\lambda'_2+\delta_1\lambda'_1-\delta_1^2).
$$

Since $\Delta_1=[\mg{1,1}^2/S_2]$, its Chow ring is a quotient of
$\CH(BG_2)=\mathbb Z[\beta_1,\beta_2,\gamma]/(2\gamma,\gamma(\beta_1+\gamma))$, where
$\gamma$ is the sign class and $\beta_1,\beta_2$ are
the first and second Chern classes of the standard representation.
By \cite[Lemma 6.1]{Larson}, we have
$$
j^*\delta_1=\gamma-\beta_1=\gamma-j^*\lambda'_1,
$$
and hence $\gamma=j^*(\delta_1+\lambda'_1)$.
To see the second equality, let $r_1=10\lambda'_1-2\delta_1$
and $r_2=2{\lambda'_1}^2-24\lambda'_2+\delta_1\lambda'_1-\delta_1^2$ (both relations on
$\mg2^{\mathrm{ct}}$). Using $j_*\gamma=j_*j^*(\delta_1+\lambda'_1)$, we have
$$
j_*\gamma-12({\lambda'_1}^2-2\lambda'_2)=\delta_1(\delta_1+\lambda'_1)-
12({\lambda'_1}^2-2\lambda'_2)=-r_2-\lambda'_1 r_1=0.
$$
\end{proof}

\begin{proposition}\label{test family pullbacks}
Pullback along $f$ is given by
$$
\begin{aligned}
\lambda_1&\mapsto \lambda'_1,&
\lambda_2&\mapsto \lambda'_2,&
\lambda_3&\mapsto 0,\\
\alpha_{12}&=(p_{12})_*1\mapsto \lambda'_2-\delta_1^2-2\lambda'_1\delta_1,&
\alpha_{12}^{(1)}&=(p_{12})_*t\mapsto \lambda'_1\delta_1^2+2\lambda'_2\delta_1,\\
\alpha_{111}&=(p_{111})_*1\mapsto \lambda'_2\delta_1,&
\alpha_{111}^{(1)}&=(p_{111})_*d_1\mapsto \lambda'_2\delta_1(\delta_1+\lambda'_1).
\end{aligned}
$$
\end{proposition}

\begin{proof}
The Hodge bundle pulls back to $\mathcal O\oplus\mathbb E_2$, which gives the formulas for the
$\lambda_i$. To pull back $\alpha_{12}=(p_{12})_*1$ and $\alpha_{12}^{(1)}=(p_{12})_*t$, consider the fiber product
$$
\begin{tikzcd}[column sep=large, row sep=large]
\mg2^{\mathrm{ct}}\times_{\mathcal A_3}(\mathcal A_1\times\mathcal A_2)
\arrow[r,"g"]\arrow[d,"q"'] & \mathcal A_1\times\mathcal A_2\arrow[d,"p_{12}"]\\
\mg2^{\mathrm{ct}}\arrow[r,"f"'] & \mathcal A_3.
\end{tikzcd}
$$
Note that we have
$$
\mg2^{\mathrm{ct}}\times_{\mathcal A_3}(\mathcal A_1\times\mathcal A_2)
=\mg2^{\mathrm{ct}}\sqcup(\mg{1,1}\times\mg{1,1}).
$$
Indeed, the fiber product parametrizes choices of an elliptic factor $E$ in $E_0\times J(C)$.
Choosing $E=E_0$ gives $\mg2^{\mathrm{ct}}$, while choosing $E=E_1$ in
$E_0\times E_1\times E_2$ gives $\mg{1,1}\times\mg{1,1}$.
The map $q$ decomposes into $q=(\id,q')$ where $q'$ is the composition
$$
q':\mg{1,1}\times\mg{1,1}\xrightarrow{\nu}[\mg{1,1}^2/S_2]=\Delta_1\xrightarrow{j}\mg2^{\mathrm{ct}}.
$$
The map $g$ decomposes into $g=(g_1,g_2)$ where
$g_1:\mg2^{\mathrm{ct}}\rightarrow\mathcal A_1\times\mathcal A_2$ is the map
$C\mapsto(E_0,J(C))$ and
$g_2:\mg{1,1}\times\mg{1,1}\rightarrow\mathcal A_1\times\mathcal A_2$ is the map
$(E_1,E_2)\mapsto(E_1,E_0\times E_2)$.

Let $L_1,L_2$ be the Hodge bundles on the first and second factors, respectively,
of $\mg{1,1}\times\mg{1,1}$. Note that $L_1$ is pulled back from
$\mathcal A_1\times\mathcal A_2$ along $g_2$, and so we continue referring to its
first Chern class as $t$. Let $\Gamma$ be the sign line bundle
on $\Delta_1=[(\mg{1,1}\times\mg{1,1})/S_2]$, with $\gamma=c_1(\Gamma)$.
By Lemma \ref{Kodaira-Spencer normal}, the normal bundle along $p_{12}$ is
$$
N_{p_{12}}=\mathbb E_1^\vee\otimes\mathbb E_2^\vee.
$$
Along $q$, the normal bundle on $\mg2^{\mathrm{ct}}$ is zero and on
$\mg{1,1}\times\mg{1,1}$ is $L_1^{-1}\otimes L_2^{-1}$. Hence
$$
g_1^*N_{p_{12}}\cong\mathbb E_2^\vee,\qquad
g_2^*N_{p_{12}}\cong L_1^{-1}\oplus(L_1^{-1}\otimes L_2^{-1}).
$$
Thus the excess bundle on $\mg2^{\mathrm{ct}}$ is $\mathbb E_2^\vee$, whose top
Chern class is $\lambda'_2$, and on $\mg{1,1}\times\mg{1,1}$ is $L_1^{-1}$, whose top
Chern class is $-t$.
Then excess intersection gives
\begin{equation}\label{f pullback}
\begin{aligned}
f^*\alpha_{12}&=f^*(p_{12})_*1=\id_*(\lambda'_2\cdot1)+q'_*(-t\cdot1)=\lambda'_2-j_*\nu_*t,\\
f^*\alpha_{12}^{(1)}&=f^*(p_{12})_*t=\id_*(\lambda'_2\cdot0)+q'_*(-t\cdot t)=-j_*\nu_*t^2,
\end{aligned}
\end{equation}
so we only need $j_*\nu_*t$ and $j_*\nu_*t^2$. \cite[Lemma 7.1]{Larson}, keeping in mind that
our $\nu$ is his $\pi$ and our $t$ is his $t_1$, gives
$$
\begin{aligned}
\nu_*t&=j^*\lambda'_1+\gamma,\\
\nu_*t^2&=(j^*\lambda'_1)(j^*\lambda'_1+\gamma)-2j^*\lambda'_2.
\end{aligned}
$$
The projection formula, using $2\gamma=0$, therefore gives
$$
\begin{aligned}
j_*\nu_*t&=\lambda'_1\delta_1+j_*\gamma=\delta_1^2+2\lambda'_1\delta_1,\\
j_*\nu_*t^2&={\lambda'_1}^2\delta_1-\lambda'_1 j_*\gamma-2\lambda'_2\delta_1=-\lambda'_1\delta_1^2-2\lambda'_2\delta_1.
\end{aligned}
$$
Substituting into \eqref{f pullback} gives the claimed formulas for $f^*\alpha_{12}$ and $f^*\alpha_{12}^{(1)}$.

Finally, consider the Cartesian diagram
\begin{equation}\label{test family square}
\begin{tikzcd}[column sep=large, row sep=large]
\Delta_1\arrow[r,"h"]\arrow[d,"j"'] & A_{1,1,1}\arrow[d,"p_{111}"]\\
\mg2^{\mathrm{ct}}\arrow[r,"f"'] & \mathcal A_3.
\end{tikzcd}
\end{equation}
By Lemma \ref{Kodaira-Spencer normal}, the excess bundle is $j^*\mathbb E_2^\vee$, with top
Chern class $j^*\lambda'_2$. Therefore
\begin{equation}\label{test family excess}
f^*(p_{111})_*z=\lambda'_2j_*h^*z\qquad(z\in\CH(A_{1,1,1})).
\end{equation}
Since $h^*d_1=\gamma$, taking $z=1,d_1$ gives the formulas for $f^*\alpha_{111}$ and $f^*\alpha_{111}^{(1)}$.
\end{proof}

\section{The integral Chow ring of $\mathcal A_3$}\label{final ring section}

We now assemble the pieces. By Lemma \ref{test family restriction lemma}, it suffices to
show that the simultaneous restriction to each of $\CH(\mathcal A_3^{\mathrm{ind}})$,
$\CH(A_{1,2}^\circ)$, $\CH(A_{1,1,1})$, $\CH(\mg2^{\mathrm{ct}})$, and
$\CH(\mathcal A_3)_{\mathbb Q}$ is injective. Once we do this, we will only need to determine
the generators of $\CH(\mathcal A_3)$ and then intersect the kernels of the induced maps to each
of the above Chow rings.

\begin{proposition}\label{simultaneous restriction lemma}
The homomorphism
\begin{equation}\label{simultaneous restriction map}
\begin{aligned}
\CH(\mathcal A_3)&\rightarrow
\CH(\mathcal A_3^{\mathrm{ind}})\times\CH(A_{1,2}^\circ)\times\CH(A_{1,1,1})
\times\CH(\mg2^{\mathrm{ct}})\times\CH(\mathcal A_3)_{\mathbb Q},\\
\omega&\mapsto
\left(\omega|_{\mathcal A_3^{\mathrm{ind}}},p_{12}^*\omega,p_{111}^*\omega,
f^*\omega,\omega_{\mathbb Q}\right)
\end{aligned}
\end{equation}
is injective.
\end{proposition}

\begin{proof}
Let $\omega$ lie in the kernel, and assume that $\omega$ is homogeneous.
By Lemma \ref{A_12 detection}, the rational vanishing of $\omega$ plus its vanishing on
$\mathcal A_3^{\mathrm{ind}}$ and $A_{1,2}^\circ$ imply that
$\omega|_{\mathcal A_3\setminus A_{1,1,1}}=0$.
Then excision says that $\omega=(p_{111})_*z$ for some $z\in\CH(A_{1,1,1})$,
and by self-intersection $ez=p_{111}^*\omega=0$. That is, $z\in\Ann(e)$.

By van der Geer's computation \cite[Propositions 3.1, 3.2 and 3.13]{vdG},
$$
\alpha_{111}=35\lambda_1^3/2\neq0
\quad\text{in }\CH(\mathcal A_3)_\mathbb Q.
$$
If $\deg z=0$, then $z=n\in\mathbb Z$, and $0=\omega_{\mathbb Q}=n\alpha_{111}$ implies $n=0$.
Now suppose $\deg z>0$. By Lemma \ref{A_1,1,1 annihilator}, we have $2z=0$. Since
$\omega=(p_{111})_*z$, this gives
\begin{equation}\label{2 omega vanishing}
2\omega=0.
\end{equation}
So to show that $\omega$ vanishes, it would suffice to also show $3\omega=0$.

Since $\Ann(e)=(36,d_1)$ (Lemma \ref{A_1,1,1 annihilator}) and $36d_1=36d_2=d_1d_2=0$
(Proposition \ref{A_1,1,1 Chow ring}), we can write
$$
z=36P(c_1,c_2,c_3)+d_1 F(d_1,c_1,c_2,c_3),
$$
where $P$ has positive degree. In the square from \eqref{test family square}, we have
$$
h^*d_1=\gamma=j^*(\delta_1+\lambda'_1),\quad
h^*c_1=j^*\lambda'_1,\quad h^*c_2=j^*\lambda'_2,\quad h^*c_3=0.
$$
Put $\sigma=j_*h^*z\in\CH(\mg2^{\mathrm{ct}})$. Then
\begin{equation}\label{alternate sigma form}
\sigma=36\delta_1P(\lambda'_1,\lambda'_2,0)+\delta_1(\delta_1+\lambda'_1)F
(\delta_1+\lambda'_1,\lambda'_1,\lambda'_2,0).
\end{equation}
By \eqref{test family excess}, we have $0=f^*\omega=\lambda'_2\sigma$.
Since $P$ has positive degree, $\sigma$ belongs to the ideal
$$
(36\delta_1\lambda'_1,36\delta_1\lambda'_2,
\delta_1(\delta_1+\lambda'_1))\subseteq\CH(\mg2^{\mathrm{ct}}).
$$
The computation in \texttt{proposition-12-1.m2} verifies that multiplication by $\lambda'_2$ is
injective on this ideal. Hence $\lambda'_2\sigma=0$ implies $\sigma=0$, and so
\begin{equation}\label{sigma vanishing}
\sigma=0.
\end{equation}

Recall the product morphism $p_{12}:\mathcal A_1\times\mathcal A_2\rightarrow\mathcal A_3$.
Under $\mathcal A_2\cong\mg2^{\mathrm{ct}}$, it is given by
$(E,C)\mapsto E\times J(C)$. Recall
$$
\pi:\mathcal A_1\times [\mathcal A_1^2/S_2]\rightarrow[\mathcal A_1^3/S_3]=A_{1,1,1}
$$
from Definition \ref{pi definition}. Note that $\Delta_1\cong[\mathcal A_1^2/S_2]$, so we will
view the domain of $\pi$ as $\mathcal A_1\times\Delta_1$. These maps fit into the Cartesian diagram
$$
\begin{tikzcd}[column sep=large, row sep=large]
\mathcal A_1\times\Delta_1\arrow[r,"\pi"]\arrow[d,"1\times j"'] & A_{1,1,1}\arrow[d,"p_{111}"]\\
\mathcal A_1\times\mg2^{\mathrm{ct}}\arrow[r,"p_{12}"'] & \mathcal A_3.
\end{tikzcd}
$$
Restricting the first factor to the fixed elliptic curve $E_0\in\mathcal A_1$ recovers the diagram
\eqref{test family square}. Let
$q:\mathcal A_1\times\mg2^{\mathrm{ct}}\rightarrow\mg2^{\mathrm{ct}}$ be the projection. Since
$\CH(\mathcal A_1\times\mg2^{\mathrm{ct}})=\CH(\mg2^{\mathrm{ct}})[t]/(12t)$, we will suppress
pullback by $q$.

We claim that
\begin{equation}\label{sigma equality}
(1\times j)_*\pi^*z=\sigma.
\end{equation}
The pullbacks of $d_1,c_1,c_2,c_3$ along $\pi$ are the restrictions of
$$
\delta_1+\lambda'_1,\qquad t+\lambda'_1,\qquad t\lambda'_1+\lambda'_2,\qquad t\lambda'_2,
$$
respectively. Applying the projection formula to
$z=36P(c_1,c_2,c_3)+d_1 F(d_1,c_1,c_2,c_3)$ therefore gives
$$
\begin{aligned}
(1\times j)_*\pi^*z&=36\delta_1P(t+\lambda'_1,t\lambda'_1+\lambda'_2,t\lambda'_2)\\
&\quad+\delta_1(\delta_1+\lambda'_1)
F(\delta_1+\lambda'_1,t+\lambda'_1,t\lambda'_1+\lambda'_2,t\lambda'_2).
\end{aligned}
$$
Notice that this is exactly our expression for $\sigma$ from \eqref{alternate sigma form}, but
with $t=0$:
$$
\sigma=36\delta_1P(\lambda'_1,\lambda'_2,0)+\delta_1(\delta_1+\lambda'_1)F
(\delta_1+\lambda'_1,\lambda'_1,\lambda'_2,0).
$$
Hence $\sigma$ differs from $(1\times j)_*\pi^*z$ by a multiple of $t$.
Now notice that the first term in each of $(1\times j)_*\pi^*z$ and $\sigma$ is
divisible by 36 and that their second terms both contain
$\delta_1(\delta_1+\lambda'_1)=12({\lambda'_1}^2-2\lambda'_2)$
(Proposition \ref{M_2^ct Chow ring prop}). Thus their difference must vanish since $12t=0$.
This proves the claim.

Thus, using $\sigma=0$ from \eqref{sigma vanishing} and $(1\times j)_*\pi^*z=\sigma$ from
\eqref{sigma equality},
$$
3\omega=(p_{111})_*3z=(p_{111})_*\pi_*\pi^*z=(p_{12})_*(1\times j)_*\pi^*z=(p_{12})_*\sigma=0.
$$
Since $2\omega=0$ from \eqref{2 omega vanishing}, we conclude that $\omega=0$.
\end{proof}

\begin{proposition}\label{A_3 generators}
The ring $\CH(\mathcal A_3)$ is generated by $\lambda_1,\lambda_2,\lambda_3,\alpha_{12},\alpha_{12}^{(1)},\alpha_{111},\alpha_{111}^{(1)}$.
\end{proposition}

\begin{proof}
The $\lambda_i$ generate $\CH(\mathcal A_3^{\mathrm{ind}})$ by Proposition \ref{A_ind Chow ring}.
On $A_{1,2}^\circ$, Proposition \ref{A_12 Chow ring} gives
$$
\lambda'_1=p_{12}^*\lambda_1-t,\qquad
\lambda'_2=p_{12}^*\alpha_{12}-p_{12}^*\lambda_1 t,\qquad
t^2=p_{12}^*(\alpha_{12}-\lambda_2).
$$
Thus $\CH(A_{1,2}^\circ)$ is generated by $1,t$ as a module over polynomials in
$\lambda_1,\lambda_2,\lambda_3,\alpha_{12}$, and so excision and the projection formula show that
$\CH(\mathcal A_3\setminus A_{1,1,1})$ is generated by
$\lambda_1,\lambda_2,\lambda_3,\alpha_{12},\alpha_{12}^{(1)}$. So we have reduced to understanding the image
of $(p_{111})_*$.

Put $\theta=\alpha_{12}-\lambda_1^2-\lambda_2$. By \eqref{i^* equalities},
$p_{111}^*\theta=d_1^2+2d_2$. Since $2d_1=3d_2=0$, we have
$$
d_1^2=3p_{111}^*\theta,\qquad d_2=2p_{111}^*\theta.
$$
Thus $\CH(A_{1,1,1})$ is generated by $1,d_1$ as a module over polynomials in
$\lambda_1,\lambda_2,\lambda_3,\alpha_{12}, \alpha_{12}^{(1)}$. The projection formula now shows that the image of
$(p_{111})_*$ is generated by $\alpha_{111}=(p_{111})_*1$ and $\alpha_{111}^{(1)}=(p_{111})_*d_1$. Hence we conclude that
$\CH(\mathcal A_3)$ is generated by $\lambda_1,\lambda_2,\lambda_3,\alpha_{12},\alpha_{12}^{(1)},\alpha_{111},\alpha_{111}^{(1)}$.
\end{proof}

\begin{definition}\label{five evaluation maps}
Put
$$
R=\mathbb Z[\lambda_1,\lambda_2,\lambda_3,\alpha_{12},\alpha_{12}^{(1)},\alpha_{111},\alpha_{111}^{(1)}].
$$
Let $\phi_{\mathrm{ind}}$, $\phi_{12}$, $\phi_{111}$, $\phi_{\mathrm{ct}}$ and
$\phi_{\mathbb Q}$ be the homomorphisms obtained by composing
$R\rightarrow\CH(\mathcal A_3)$ with the projection to each term in
\eqref{simultaneous restriction map}. Write $K_{\mathrm{ind}}$, $K_{12}$, $K_{111}$,
$K_{\mathrm{ct}}$ and $K_{\mathbb Q}$ for their kernels in $R$.

The map $\phi_{\mathrm{ind}}$ sends $\alpha_{12},\alpha_{12}^{(1)},\alpha_{111},\alpha_{111}^{(1)}$ to zero and sends the
$\lambda_i$ to the $\lambda_i$.
By Proposition \ref{A_12 Chow ring}, $\phi_{12}$ is given by
$$
\lambda_1\mapsto t+\lambda'_1,\quad \lambda_2\mapsto t\lambda'_1+\lambda'_2,\quad
\lambda_3\mapsto t\lambda'_2,\quad \alpha_{12}\mapsto e,\quad \alpha_{12}^{(1)}\mapsto et,\quad\alpha_{111}\mapsto 0,
\quad\alpha_{111}^{(1)}\mapsto0.
$$
By Proposition \ref{A_1,1,1 restrictions} and Lemma \ref{A_1,1,1 annihilator}, the map
$\phi_{111}$ sends the $\lambda_i$ to the $c_i$, and on
$\alpha_{12},\alpha_{12}^{(1)},\alpha_{111}, \alpha_{111}^{(1)}$ is given by
$$
\begin{aligned}
\alpha_{12}&\mapsto c_1^2+c_2+d_1^2+2d_2,\\
\alpha_{12}^{(1)}&\mapsto c_1^3-2c_1c_2+3c_3+d_1^2c_1,\\
\alpha_{111}&\mapsto c_3-c_1c_2+d_1(c_1^2+c_2)+d_1^3,\\
\alpha_{111}^{(1)}&\mapsto d_1e=0.
\end{aligned}
$$
By Proposition \ref{test family pullbacks}, the map $\phi_{\mathrm{ct}}$ is
the following:
$$
\begin{aligned}
\lambda_1&\mapsto \lambda'_1,&\lambda_2&\mapsto \lambda'_2,&\lambda_3&\mapsto 0,\\
\alpha_{12}&\mapsto \lambda'_2-\delta_1^2-2\lambda'_1\delta_1,&
\alpha_{12}^{(1)}&\mapsto \lambda'_1\delta_1^2+2\lambda'_2\delta_1,&
\alpha_{111}&\mapsto \lambda'_2\delta_1,\\
\alpha_{111}^{(1)}&\mapsto \lambda'_2\delta_1(\delta_1+\lambda'_1).
\end{aligned}
$$
Finally, identifying $\CH(\mathcal A_3)_{\mathbb Q}$ with
$\mathbb Q[\lambda_1]/(\lambda_1^4)$ by \cite[Propositions 3.1, 3.2, 3.13 and Corollary 3.4]{vdG},
we have
$$
\begin{aligned}
\phi_{\mathbb Q}(\lambda_1)&=\lambda_1,&
\phi_{\mathbb Q}(\lambda_2)&=\lambda_1^2/2,&
\phi_{\mathbb Q}(\lambda_3)&=0,\\
\phi_{\mathbb Q}(\alpha_{12})&=21\lambda_1^2/2,&
\phi_{\mathbb Q}(\alpha_{12}^{(1)})&=0,&
\phi_{\mathbb Q}(\alpha_{111})&=35\lambda_1^3/2,&
\phi_{\mathbb Q}(\alpha_{111}^{(1)})&=0.
\end{aligned}
$$
\end{definition}

\begin{lemma}\label{rational evaluation kernel}
The kernel of $\phi_{\mathbb Q}:R\rightarrow\CH(\mathcal A_3)_{\mathbb Q}$ is
$$
K_{\mathbb Q}=\left(\lambda_1^2-2\lambda_2,\lambda_2^2,\lambda_3,
\alpha_{12}-21\lambda_2,\alpha_{12}^{(1)},\alpha_{111}-35\lambda_1\lambda_2,\alpha_{111}^{(1)}\right).
$$
\end{lemma}

\begin{proof}
Let $I$ be the ideal generated by these elements. Each generator of $I$ maps to zero,
and every class modulo $I$ is a $\mathbb Z$-linear combination of
$$
1,\qquad\lambda_1,\qquad\lambda_2,\qquad\lambda_1\lambda_2.
$$
Their images are, respectively,
$$
1,\qquad\lambda_1,\qquad\frac{\lambda_1^2}{2},\qquad
\frac{\lambda_1^3}{2},
$$
which are linearly independent over $\mathbb Q$ in
$\CH(\mathcal A_3)_{\mathbb Q}=\mathbb Q[\lambda_1]/(\lambda_1^4)$. Hence the induced
map from the quotient by $I$ is injective, proving the claimed equality.
\end{proof}

\begin{theorem}\label{A_3 Chow ring}
We have
$$
\CH(\mathcal A_3)=\frac{\mathbb Z[\lambda_1,\lambda_2,\lambda_3,\alpha_{12},\alpha_{12}^{(1)},
\alpha_{111}]}{I_{16}},
$$
where
$\deg(\lambda_1,\lambda_2,\lambda_3,\alpha_{12},\alpha_{12}^{(1)},\alpha_{111})=(1,2,3,2,3,3)$
and the relations in $I_{16}$ are listed in Appendix \ref{A_3 relations appendix}. In particular,
the order of $\lambda_3$ is 504.
\end{theorem}

\begin{proof}
We use the criterion of Lemma \ref{test family restriction lemma}.
By Proposition \ref{A_3 generators}, the natural homomorphism
$R\rightarrow\CH(\mathcal A_3)$ is surjective.
Since the map of Proposition \ref{simultaneous restriction lemma} is injective, we see that
$$
\ker\left(R\rightarrow\CH(\mathcal A_3)\right)
=K_{\mathrm{ind}}\cap K_{12}\cap K_{111}\cap K_{\mathrm{ct}}\cap K_{\mathbb Q}.
$$
The computation in \texttt{theorem-12-5.m2} then gives
$$
K_{\mathrm{ind}}\cap K_{12}\cap K_{111}\cap K_{\mathrm{ct}}\cap K_{\mathbb Q}
=I_{16}R+\left(\alpha_{111}^{(1)}-(5562\lambda_1^2+30\lambda_2-54\alpha_{12})\alpha_{12}
-972\lambda_1\alpha_{12}^{(1)}\right),
$$
where $I_{16}$ is generated by the 16 polynomials of Appendix \ref{A_3 relations appendix},
and also computes the order of $\lambda_3$. Eliminating $\alpha_{111}^{(1)}$ gives the claimed presentation.
\end{proof}

\begin{corollary}\label{Torelli corollary}
Pullback along the Torelli morphism is surjective and identifies
$$
\CH(\mg3)=\CH(\mathcal A_3)/(\alpha_{12}, \alpha_{12}^{(1)}, \alpha_{111},252\lambda_1^2).
$$
\end{corollary}

\begin{proof}
Observe that all terms in all relations on $\mathcal A_3$ contain at least one of
$\alpha_{12},\alpha_{12}^{(1)},\alpha_{111}$, with the exception of \eqref{app:A3:1},
\eqref{app:A3:2}, \eqref{app:A3:5}, \eqref{app:A3:6}, \eqref{app:A3:8}.
\texttt{corollary-12-6.m2} verifies that the remaining relations, with
$\alpha_{12}=\alpha_{12}^{(1)}=\alpha_{111}=0$, generate the same ideal as the relations on
$\mathcal A_3^{\mathrm{ind}}$ from Proposition \ref{A_ind Chow ring}. Hence the kernel of
restriction to $\mathcal A_3^{\mathrm{ind}}$ is $(\alpha_{12},\alpha_{12}^{(1)},\alpha_{111})$,
and so combining with Theorem \ref{M_3 Chow ring} proves the claim.
\end{proof}

\appendix

\section{Guide to the computational supplement}\label{computational supplement appendix}

The computational supplement consists of the fifteen Macaulay2 scripts listed below.
Instructions are provided in \texttt{README.md}.

\begin{center}
\begin{tabular}{ll}
\hline
File & Result \\
\hline
\texttt{lemma-3-9.m2} & Lemma \ref{ordinary and special Chow rings} \\
\texttt{lemma-5-2.m2} & Lemma \ref{112 annihilator} \\
\texttt{lemma-5-3.m2} & Lemma \ref{excess chern classes} \\
\texttt{proposition-5-5.m2} & Proposition \ref{Y_1 and Y_13 pushforwards} \\
\texttt{proposition-5-6.m2} & Proposition \ref{Y_2 pushforward} \\
\texttt{proposition-5-7.m2} & Proposition \ref{no relations} \\
\texttt{proposition-7-6.m2} & Proposition \ref{A_ind Chow ring} \\
\texttt{theorem-7-8.m2} & Theorem \ref{M_3 Chow ring} \\
\texttt{theorem-8-3.m2} & Theorem \ref{A_12 annihilator} \\
\texttt{proposition-9-5.m2} & Proposition \ref{BG Chow ring} \\
\texttt{proposition-9-7.m2} & Proposition \ref{A_1,1,1 Chow ring} \\
\texttt{lemma-10-3.m2} & Lemma \ref{A_1,1,1 annihilator} \\
\texttt{proposition-12-1.m2} & Proposition \ref{simultaneous restriction lemma} \\
\texttt{theorem-12-5.m2} & Theorem \ref{A_3 Chow ring} \\
\texttt{corollary-12-6.m2} & Corollary \ref{Torelli corollary} \\
\hline
\end{tabular}
\end{center}

\section{Relations for $\mathcal A_3$}\label{A_3 relations appendix}

The presentation
$$
\CH(\mathcal A_3)=\mathbb Z[\lambda_1,\lambda_2,\lambda_3,\alpha_{12},\alpha_{12}^{(1)},\alpha_{111}]/I_{16}
$$
has the following 16 relations, listed in increasing degree:
\begingroup
\raggedbottom
\allowdisplaybreaks[4]
\setlength{\jot}{2pt}
\begin{gather}
48\alpha_{12}-504\lambda_1^2=0,\label{app:A3:1}\\[10pt]
312\alpha_{12}-3264\lambda_1^2-24\lambda_2=0,\label{app:A3:2}\\[10pt]
-10\lambda_1\alpha_{12}-2\alpha_{12}^{(1)}+6\alpha_{111}=0,\label{app:A3:3}\\[10pt]
90\lambda_1\alpha_{12}+30\alpha_{12}^{(1)}-54\alpha_{111}=0,\label{app:A3:4}\\[10pt]
\begin{aligned}
&15372\lambda_1^3+224\lambda_1\lambda_2+28\lambda_3\\
&\quad +1062\lambda_1\alpha_{12}+798\alpha_{12}^{(1)}-1522\alpha_{111}=0,
\end{aligned}\label{app:A3:5}\\[10pt]
44576\lambda_1^3+448\lambda_1\lambda_2-224\lambda_3
+1024\alpha_{12}^{(1)}-2560\alpha_{111}=0,\label{app:A3:6}\\[10pt]
-2(\alpha_{12}-\lambda_2-180\lambda_1^2)\alpha_{12}
+62\lambda_1\alpha_{12}^{(1)}-6\lambda_1\alpha_{111}=0,\label{app:A3:7}\\[10pt]
\begin{aligned}
&-18\lambda_1\lambda_3-162\lambda_1^2\lambda_2-14580\lambda_1^4\\
&\quad -(21330\lambda_1^2+132\lambda_2-36\alpha_{12})\alpha_{12}-4158\lambda_1\alpha_{12}^{(1)}=0,
\end{aligned}\label{app:A3:8}\\[10pt]
12\lambda_1\alpha_{111}=0,\label{app:A3:9}\\[10pt]
\begin{aligned}
&(12150\lambda_1^3+63\lambda_1\lambda_2-\lambda_3
-63\lambda_1\alpha_{12})\alpha_{12}\\
&\quad +(2070\lambda_1^2+3\alpha_{12}-2\lambda_2)\alpha_{12}^{(1)}-2\alpha_{12}\alpha_{111}=0,
\end{aligned}\label{app:A3:10}\\[10pt]
24\lambda_2\alpha_{111}=0,\label{app:A3:11}\\[10pt]
3(\alpha_{12}-\lambda_1^2-\lambda_2)\alpha_{111}=0,\label{app:A3:12}\\[10pt]
\begin{aligned}
&(2070\lambda_1^2+3\alpha_{12}-2\lambda_2)
(\alpha_{12}-\lambda_2-180\lambda_1^2)\alpha_{12}\\
&\quad +(-62370\lambda_1^3+135\lambda_1\lambda_2-\lambda_3
-171\lambda_1\alpha_{12})\alpha_{12}^{(1)}\\
&\quad -(8\lambda_1^3+11\lambda_1\lambda_2+33\lambda_3)\alpha_{111}\\
&\quad -\lambda_2\left((5562\lambda_1^2+30\lambda_2-54\alpha_{12})\alpha_{12}
+972\lambda_1\alpha_{12}^{(1)}\right)=0,
\end{aligned}\label{app:A3:13}\\[10pt]
\begin{aligned}
&{\alpha_{12}^{(1)}}^2+36\lambda_1\alpha_{12} \alpha_{12}^{(1)}
-(\alpha_{12}-\lambda_2-180\lambda_1^2)\alpha_{12}^2\\
&\quad -(6\lambda_1^3+9\lambda_1\lambda_2+9\lambda_3)\alpha_{111}\\
&\quad -\lambda_2\left((5562\lambda_1^2+30\lambda_2-54\alpha_{12})\alpha_{12}
+972\lambda_1\alpha_{12}^{(1)}\right)=0,
\end{aligned}\label{app:A3:14}\\[10pt]
(\alpha_{12}^{(1)}-\lambda_1\alpha_{12}-3\lambda_3+3\lambda_1\lambda_2)\alpha_{111}
=0,\label{app:A3:15}\\[10pt]
\begin{aligned}
&\alpha_{111}^2-(\lambda_3-\lambda_1\lambda_2)\alpha_{111}\\
&\quad -\lambda_2\left((5562\lambda_1^2+30\lambda_2-54\alpha_{12})\alpha_{12}
+972\lambda_1\alpha_{12}^{(1)}\right)=0.
\end{aligned}\label{app:A3:16}
\end{gather}
\endgroup

\newpage
\begingroup\small
\newcommand{\etalchar}[1]{$^{#1}$}

\endgroup

\end{document}